\documentclass[a4paper,12pt]{article}

\usepackage[english]{babel}
\usepackage{amsmath}
\usepackage{amssymb}
\usepackage{epsfig}
\usepackage{theorem}
\usepackage{enumerate}
\usepackage{hyperref}
\usepackage{calc}
\usepackage{amsfonts}
\usepackage{ifthen}

\newcounter{count}[section]
\renewcommand{\thecount}{\thesection.\arabic{count}}

\newcommand{\equa}[1]{\addtocounter{count}{1}\begin{equation}#1\end{equation}}
\newcommand{\eqna}[1]{\addtocounter{count}{1}\begin{eqnarray}#1\end{eqnarray}}

\newcounter{tictac1}
\newcounter{tictac2}
\newcounter{tictac3}
\newcounter{tictac4}

\newenvironment{fleuvea}{
   \begin{list}{\rm{$\textbf{(\alph{tictac2})}$} }{\usecounter{tictac2}
\leftmargin 1cm\labelwidth 2em}}{\end{list}}

\newenvironment{fleuvei}{
   \begin{list}{\rm{$\textbf{(\roman{tictac4})}$} }{\usecounter{tictac4}
\leftmargin 1cm\labelwidth 2em}}{\end{list}}

\newlength{\Longueur}
\newenvironment{remas1}[1]{\noindent\textbf{Remarks}\vspace{-.25cm}
\ifthenelse{#1=0}{\settowidth{\Longueur}{00.00}}{\addtocounter{count}{#1}
\settowidth{\Longueur}{\thecount}\addtocounter{count}{-#1}}
\begin{list}{\addtocounter{count}{1}{\bf
      \thecount}}{\setlength{\labelsep}{.2cm}\setlength{\labelwidth}
{\the\Longueur}\setlength{\leftmargin}{\labelwidth+\labelsep}}}{\end{list}}

\theoremstyle{break}
\newtheorem{theo}[count]{Theorem}
\newtheorem{coro}[count]{Corollary}
\newtheorem{prop}[count]{Proposition}

\newtheorem{defi}[count]{Definition}
\newtheorem{lemm}[count]{Lemma}
\theorembodyfont{\rmfamily}
\newtheorem{rema}[count]{Remark}

\newcommand{\var}[1]{{\rm var}_{#1}}

\newcommand{\reff}[1]{(\ref{#1})}

\newcommand{\lat}{\mathbb{Z}}

\def\1{\rlap{\mbox{\small\rm 1}}\kern.15em 1}
\def\ind#1{\1_{#1}}
\def\build#1_#2^#3{\mathrel{\mathop{\kern 0pt#1}\limits_{#2}^{#3}}}
\def\tend#1#2#3{\build\hbox to 12mm{\rightarrowfill}_{#1\rightarrow
#2}^{#3}}

{}
\def\converge#1#2#3{\build\hbox to
15mm{\rightarrowfill}_{\hbox{\scriptsize #3}}^{#1\rightarrow #2}}
\def\converg#1#2#3{\build\hbox to
15mm{\rightarrowfill}_{\hbox{\scriptsize #3}}^{#1\uparrow #2}}
\def\embf#1{\emph{\bf #1}}

\title{\bf{Chains with complete connections: General theory,
uniqueness, loss of memory and mixing properties}}
\author{Roberto Fern\'andez\thanks{roberto.fernandez@univ-rouen.fr}\\
Gr\'egory Maillard\thanks{gregory.maillard@univ-rouen.fr}\\
\small{{\it Laboratoire de Math\'ematiques Rapha{\"e}l Salem}}\\ \small{\it{UMR 6085 CNRS-Universit\'e de Rouen}}\\
\small{\it{Site Colbert F-76821 Mont Saint Aignan, France}}}
\begin{document}
\maketitle {\begin{abstract} We introduce an statistical mechanical formalism for the study of discrete-time stochastic
processes with which we prove:  (i) General properties of extremal chains, including triviality on the tail
$\sigma$-algebra, short-range correlations, realization via infinite-volume limits and ergodicity.  (ii) Two new
sufficient conditions for the uniqueness of the consistent chain.  The first one is a transcription of a criterion due
to Georgii for one-dimensional Gibbs measures, and the second one corresponds to Dobrushin criterion in statistical
mechanics.  (iii) Results on loss of memory and mixing properties for chains in the Dobrushin regime. These results are
complementary of those existing in the literature, and generalize the Markovian results based on the Dobrushin ergodic
coefficient.
\end{abstract}}
\thispagestyle{empty}
\section{Introduction}

\emph{Chains with complete connections} is the name coined by Onicescu and Mihoc (1935\nocite{onimih35}) for
discrete-time stochastic processes whose dependence on the past is not necessarily Markovian. The theory of these
processes has many points in common with the theory of Gibbs measures in statistical mechanics ---particularly, the
existence of phase transitions.  Nevertheless there is a clear difference, at the formal level, between both theories.
Indeed, processes are described in terms of \emph{single-site} transition probabilities, while Gibbs measures are
characterized by their conditional probabilities for \emph{arbitrary} finite regions (specifications).  In this paper
we propose a natural way to reduce this asymmetry, by introducing a statistical-mechanical framework for the study of
processes.  This framework establishes a more direct relation between both theories, which allows us to reproduce, for
chains with complete connections, a number of benchmark Gibbsian results.

We present three types of results. First, we obtain general properties of extremal chains for any type of alphabet,
namely triviality on the tail $\sigma$-algebra, short-range correlations, realization via infinite-volume limits and
ergodicity.  Second, we produce some new sufficient conditions for the uniqueness of the consistent chain.  On the one
hand, we obtain a transcription of a criterion given by Georgii (1974)\nocite{geo74} for one-dimensional Gibbs fields.
This criterion is known to be optimal for the latter, in the sense that it pinpoints the absence of phase transition
for two-body spin models with a $1/r^{2+\varepsilon}$-interaction, for all $\varepsilon>0$. The criterion imposes no
restriction on the type of alphabet.  On the other hand we prove a ``one-sided'' Dobrushin criterion, which corresponds
to a well known uniqueness criterion in statistical mechanics (see, for instance, Simon, 1993\nocite{sim93}, Chapter
V). This criterion is valid for systems with a compact metric alphabet. We exhibit simple examples where Dobrushin
criterion applies but that fall outside the scope of most other known uniqueness criteria (Harris, 1955\nocite{har55};
Iosifescu and Spataru, 1973\nocite{iosspa73}; Walters, 1955\nocite{wal75}; Berbee, 1987\nocite{ber87}; Stenflo,
2002\nocite{ste02}; Johansson and \"Oberg, 2002\nocite{johobe02}).

Our third type of results refer to loss of memory and mixing properties of chains in the Dobrushin regime.  Our
results, obtained along the lines of a similar Gibbsian theory (again we refer the reader to Chapter V of Simon, 1993),
are complementary, both in their precision and in their range of applicability, to similar results available in the
literature (Iosifescu, 1992\nocite{ios92}; Bressaud, Fern\'andez and Galves, 1999\nocite{bfg2} and references therein).
The results depend on a \emph{sensitivity matrix} that generalizes the Dobrushin ergodic coefficient of Markov chains.

Our approach is based on a notion analogous to the specifications in statistical mechanics, which we call \emph{left
interval-specifications} (LIS).  These are kernels for regions in the form of intervals which depend on the preceding
history of the process.  In contrast, Gibbsian specifications involve arbitrary finite regions and depend of the
configuration on the whole exterior of the region.  This amounts, in one dimension, to a dependence on both past and
future.  The difference is, of course, a consequence of the ``one-sidedness'' associated to a stochastic (time)
evolution, as compared with the lack of favored direction in the spatial description provides by a Gibbs measure.

The description in terms of LIS is totally equivalent to the traditional description in terms of transition
probabilities (=LIS singletons).  We show this in our first theorem.  But, as this paper illustrates, our approach has
the advantage of allowing us to ``import'', in a natural manner, notions, techniques and arguments from statistical
mechanics.  It may also be useful in the opposite direction, namely to explore the consequences of known properties of
chains for the theory of Gibbs measures.  As a step in this direction, in a companion paper (Fern\'andez and Maillard,
2003\nocite{fm03a}) we study conditions under which chains and Gibbs measures can be identified.  On a more conceptual
level, we believe that our statistical mechanical approach is more appropriate to study the general situation where
several different chains are consistent with the same transition probabilities (Bramson and Kalikow
(1993)\nocite{brakal93}, or Lacroix (2000)\nocite{lac00}). Statistical mechanics is the framework developed, precisely,
to study this phenomenon which corresponds to the appearance of (first-order) phase transitions.

\section{Preliminaries}
We consider a measurable space $(E,\mathcal{E})$ and a subset
$\Omega\subset E^{\lat}$. The exponent $\lat$ stands, in fact, for any
countable set with a total order.  The group structure of $\mathbb{Z}$
will play no role, except in Theorem \ref{th4} where
$\mathbb{Z}$ acts by isomorphisms.  The elements of $\lat$ are
called \emph{sites}, and those of $\Omega$ (\emph{admissible})
\emph{configurations}.  The space $E$ is sometimes called
\emph{alphabet}.  We endow $\Omega$ with the projection
$\mathcal{F}$ of the product $\sigma$-algebra associated to $E^\lat$.
When we invoke topological notions (e.g.\ compactness) the
$\sigma$-algebra $\mathcal{E}$ is assumed to be Borelian.  We adopt
the following notation
\begin{itemize}
\item Let $\Lambda\subset\lat$. For a configuration $\sigma\in E^\lat$
we denote $\sigma_\Lambda = (\sigma_i)_{i\in\Lambda}\in E^{\Lambda}$.
The set of admissible configurations in $\Lambda$ is $\Omega_\Lambda
\triangleq \bigl\{\sigma_\Lambda\in E^{\Lambda}: \exists\,
\omega\in \Omega \hbox{ with } \omega_\Lambda=\sigma_\Lambda \bigr\}$,
while $\mathcal{F}_{\Lambda}$ is the sub-$\sigma$-algebra of
$\mathcal{F}$ generated by the cylinders with base in
$\Omega_\Lambda$.  If $\Delta\subset\lat$ with
$\Lambda\cap\Delta=\emptyset$, $\omega_\Lambda\, \sigma_\Delta$
denotes the configuration on $\Lambda\cup\Delta$ coinciding with
$\omega_i$ for $i\in\Lambda$ and with $\sigma_i$ for $i\in\Delta$.
\item We denote $\mathcal{S}_{b}$ the set of finite intervals of
$\mathbb{Z}$.  When $\Lambda=[k,n]\in \mathcal{S}_{b}$ we shall also
use the ``sequence'' notation:
$\omega_{k}^{n}\triangleq\omega_{[k,n]}=\omega_{k}, \ldots ,
\omega_{n}$; $\Omega_k^n\triangleq\Omega_{[k,n]}$; etc.  If
$\Lambda=[k,+\infty[$, the notation will be analogous but with
$+\infty$ as upper limit.
\item If $n\in\mathbb{Z}$, $\mathcal{F}_{\leq
  n}\triangleq\mathcal{F}_{]-\infty,n]}$.  For
every $\Lambda \in \mathcal{S}_{b}$ we denote $l_\Lambda\triangleq
  \min\Lambda$; $m_\Lambda\triangleq \max\Lambda$; $\Lambda_{-} =
]-\infty,l_\Lambda - 1]$.
\item For kernels associated to a LIS (defined below), $\lim_{\Lambda
\uparrow V} f_{\Lambda}$ is the limit of the net $\left\{f_{\Lambda},
\{\Lambda\}_{\Lambda \in \mathcal{S}_{b}}, \; \Lambda \subset V ,
\subset \right\}$, for $V$ an infinite interval of $\mathbb{Z}$. If
$\mu$ a measure on $(\Omega,\mathcal{F})$ and $h$ a
$\mathcal{F}$-measurable function, we will write $\mu(h)$ instead of
$E_\mu(h)$.
\end{itemize}
\medskip

\begin{defi}[LIS]\label{lis1}
A \embf{left interval-specification} $f$ on $(\Omega, \mathcal{F})$ is
a family of probability kernels $\left\{f_{\Lambda} \right\}_{ \Lambda
\in \mathcal{S}_{b}}$, $f_{\Lambda} : \mathcal{F}_{\leq m_\Lambda}
\times \Omega \longrightarrow [0,1]$ such that for all $\Lambda$ in
$\mathcal{S}_{b}$,
\begin{fleuvea}
\item For each $\displaystyle{A \in \mathcal{F}_{\leq m_\Lambda}, \;
  f_{\Lambda}(A \mid \cdot \, )}$ is
  $\mathcal{F}_{\Lambda_{-}}$-measurable.
\item For each $\displaystyle{B \in \mathcal{F}_{\Lambda_{-}} \text{
and } \omega \in \Omega, \; f_{\Lambda}(B \mid \omega) =
\ind{B}(\omega).}$
\item For each $\displaystyle{\Delta \in \mathcal{S}_{b} : \Delta
  \supset \Lambda,}$
\equa{\label{lis3} f_{\Delta} f_{\Lambda} \;=\; f_{\Delta} \quad
  \text{on } \mathcal{F}_{\leq m_\Lambda}\;,
}
that is, $(f_{\Delta} f_{\Lambda})(h\mid\omega) =
f_{\Delta}(h\mid\omega)$ for each $\mathcal{F}_{\leq
  m_\Lambda}$-measurable function $h$ and configuration $\omega \in
\Omega$.
\end{fleuvea}
\end{defi}

These conditions are analogous to those defining a \emph{specification} in the theory of Gibbs measures (see Georgii,
1988\nocite{geo88}, for instance).  Two important differences should be highlighted, however, both being a consequence
of the ``directional'' character of the notion of process.  First, the LIS kernels act only on functions measurables
towards the left, while Gibbsian specifications have no similar constraint.  As a consequence, LIS kernels involve only
conditioning with respect to the past [property (b)], while Gibbsian kernels condition with respect to the whole
exterior of $\Lambda$.  Second, LIS kernels are defined only for intervals while Gibbsian kernels are defined for all
finite sets of sites.

Property c) is usually labeled \emph{consistency}.  There and in the sequel we adopt the standard notation for a
composition of probability kernels or of a probability kernel with a measure.  Explicitly, \reff{lis3} means that
\[
\iint h(\xi)\, f_{\Lambda}(d \xi \mid \sigma)\, f_{\Delta}(d \sigma \mid
\omega)\;=\; \int h(\sigma)\, f_{\Delta}(d \sigma \mid \omega)
\]
for each $\mathcal{F}_{\leq m_\Lambda}$-measurable function $h$ and
configuration $\omega \in \Omega$.

\begin{defi}[left interval-consistency]
A probability measure $\mu$ on $(\Omega,\; \mathcal{F} )$ is said to
be \embf{consistent} with a LIS $f$ if for each $\Lambda \in
\mathcal{S}_{b}$
\equa{
\mu f_{\Lambda}\; \;= \mu \quad \text{ on }
\mathcal{F}_{\leq m_\Lambda}\;.
}
Such a measure $\mu$ is called a \embf{chain with complete
connections}, or simply a \embf{chain}, consistent with the LIS
$f$. The family of these measures will be denoted $\mathcal{G}(f)$.
\end{defi}

\begin{remas1}{2}
\item A \emph{Markov LIS of range $k$} is a LIS such that each
function $f_{\Lambda}(A \mid \cdot \,)$ is measurable with respect to
$\mathcal{F}_{[l_\Lambda-k, l_\Lambda-1]}$, for each $A \in
\mathcal{F}_{\Lambda}$.  A chain consistent with such a LIS is a
\emph{Markov chain of range $k$}.
\item \emph{Chains with complete connections} is the original
nomenclature introduced by Onicescu and Mihoc (1935) \nocite{onimih35}. These objects have been later reintroduced
under a panoply of names, some associated to particular additional properties, others to notions later proven to be
equivalent.  Among them we mention: \emph{chains of infinite order} (Harris, 1955\nocite{har55}), \emph{$g$-measures}
(Keane, 1972\nocite{kea72}), list processes (Lalley, 1986\nocite{lal86}), \emph{uniform martingales} or \emph{random
Markov processes} (Kalikow, 1990\nocite{kal90}).
\end{remas1}

\section{Results on general framework}
We start by making the connection with the traditional definition of
chains based on singleton kernels.

%
\begin{theo}[Singleton consistency for chains]\label{listh1}
Let $\left( f_{i}\right)_{i \in \lat}$ be a family of
probability kernels $f_{i} : \mathcal{F}_{\leq i} \times \Omega
\rightarrow [0,1]$ such that for each $i \in \lat$
\begin{fleuvea}
\item For each $A \in \mathcal{F}_{\leq i}, \; f_{i}\left(A \mid \cdot
\, \right)$ is $\mathcal{F}_{\leq i-1}$-measurable.
\item For each $B \in \mathcal{F}_{\leq i-1} \text{ and } \omega \in
\Omega, \; f_{i}\left(B \mid \omega\right) = \ind{B}(\omega)$.
\end{fleuvea}
Then the LIS $f = \left\{f_{\Lambda} \right\}_{\Lambda \in
\mathcal{S}_{b}}$ defined by
\equa{\label{liss9}
f_{\Lambda} \;=\; f_{l_{\Lambda}} \, f_{l_{\Lambda}+1}
\, \cdots \, f_{m_{\Lambda}}
}
is the unique LIS such that $f_{\{i\}}=f_{i}$ for all $i\in\lat$.
Furthermore,
\equa{\label{liss7}
\mathcal{G}(f) \;=\; \Bigl\{ \mu : \mu f_{i}= \mu,
\text{ for all } i \text{ in } \lat \Bigr\}\;.
}
\end{theo}

In particular, the theorem shows that any LIS $f$ enjoys the
factorization property
\equa{\label{eq:facr}
f_{\Lambda} \;=\;
f_{\{l_{\Lambda}\}} \, f_{\{l_{\Lambda}+1\}} \, \cdots \,
f_{\{m_{\Lambda}\}}
}
on $\mathcal{F}_{\le m_{\Lambda}}$ for each
$\Lambda\in\mathcal{S}_b$. By recurrence this yields
\equa{\label{eq:facr2}
f_{[l,m]} \;=\; f_{[l,n]} \, f_{[n+1,m]}
}
for any $l,n,m\in\lat$ with $l\le n< m$.
\medskip

The following three theorems establish relations among extremality, triviality, mixing properties and infinite-volume
limits similar to those valid for Gibbs measures or, more generally, for measures consistent with specifications.
Their proofs, presented in Section \ref{s.pgf}, are patterned on the Gibbsian proofs, taking care of the one-sided
measurability of the LIS kernels.

\begin{theo}[Extremality and triviality]\label{th2}
Let $f=\left( f_{\Lambda} \right)_{\Lambda \in \mathcal{S}_{b}}$ be a
left interval-specification on $(\Omega, \mathcal{F})$. Denote by
$\mathcal{F}_{-\infty} \triangleq \bigcap_{k \in \lat}
\mathcal{F}_{\leq k}$ the \embf{tail $\sigma$-algebra}. Then
\begin{fleuvea}
\item $\mathcal{G}(f)$ is a convex set.
\item A measure $\mu$ is extreme in $\mathcal{G}(f)$ if and only if
$\mu$ is trivial on $\mathcal{F}_{-\infty}$.
\item Let $\mu \in \mathcal{G}(f)$ and $\nu \in \mathcal{P}(\Omega,
\mathcal{F})$ such that $\nu \ll \mu$.  Then
$\nu \in \mathcal{G}(f)$ if and only if there exists a
$\mathcal{F}_{-\infty}$-measurable function $h \geq 0$ such that
$\nu=h \mu$.
\item Each $\mu \in \mathcal{G}(f)$ is uniquely determined (within
$\mathcal{G}(f)$) by its restriction to the tail $\sigma$-algebra
$\mathcal{F}_{-\infty}.$
\item Two distinct extreme elements $\mu, \nu $ of $\mathcal{G}(f)$
are mutually singular on $\mathcal{F}_{-\infty}.$
\end{fleuvea}
\end{theo}

\begin{theo}[Triviality and short-range correlations]\label{theo2}
For each probability measure on $(\Omega, \mathcal{F}),$ the following
statements are equivalent.
\begin{fleuvea}
\item $\mu$ is trivial on $\mathcal{F}_{-\infty}$.
\item $\displaystyle{\lim_{\Lambda \uparrow \lat } \sup_{B \in
\mathcal{F}_{\Lambda_{-}}} \mid \mu(A \cap B) - \mu(A) \mu(B) \mid =
0}$, for all cylinder sets $A$ in $\mathcal{F}$.
\item $\displaystyle{\lim_{\Lambda \uparrow \lat } \sup_{B \in
\mathcal{F}_{\Lambda_{-}}} \mid \mu(A \cap B) - \mu(A) \mu(B) \mid =
0}$, for all $A \in \mathcal{F}$.
\end{fleuvea}
\end{theo}

\begin{theo}[Infinite volume limits]\label{th3}
Let $f$ be a LIS, $\mu$ an extreme point of $\mathcal{G}(f)$ and
$\left( \Lambda_{n} \right)_{n \geq 1}$ a sequence of regions in
$\mathcal{S}_{b}$ such that $\Lambda_n\uparrow \lat$.  Then
\begin{fleuvea}
\item $f_{\Lambda_{n}}h \rightarrow \mu(h)$ $\mu$-a.s.\ for each
bounded local function $h$ on $\Omega$
\item If $\Omega$ is a compact metric space, then for $\mu$-almost all
$\omega \in \Omega$, $f_{\Lambda_{n}}h \rightarrow \mu(h)$ for all
continuous local functions $h$ on $\Omega$.
\end{fleuvea}
\end{theo}

The following theorem is the only result in the paper where we consider translation invariance. We briefly recall the
relevant notions.  We consider the \embf{(right) shift} $\tau(i)=i+1$. (More generally, the same theory applies to any
action of $\mathbb{Z}$ on $\lat$ by isomorphisms.  In the case of $k$-shifts such theory leads to $k$-periodic
objects). The shift induces actions on configurations, measurable sets, measurable functions and measures that we
denote with the same symbol: for $\omega\in\Omega$ $\tau(\omega)=\left(\omega_{i-1}\right)_{i \in \mathbb{Z}}$, for
$A\in\mathcal{F}$, $\tau A =\{\tau\omega:\omega\in A\}$; for $h$ $\mathcal{F}$-measurable, $(\tau
h)(\omega)=h(\tau^{-1}\omega)$, and for a measure $\mu$ on $(\Omega,\mathcal{F})$, $(\tau\mu)(h)=\mu(\tau^{-1}h)$.
Objects invariant under the action of the shift are called \embf{shift-invariant}.  We denote $\mathcal{I}$ the
$\sigma$-algebra of all shift-invariant measurable sets, and $\mathcal{P}_{\text{inv}}(\Omega,\mathcal{F})$ the set of
shift-invariant probability measures on $(\Omega, \mathcal{F})$.  A measure in
$\mathcal{P}_{\text{inv}}(\Omega,\mathcal{F})$ is \embf{ergodic} if it is trivial on $\mathcal{I}$.

For $k\in\mathbb{Z}$ and $\Lambda\subset\lat$ we denote $\Lambda+k =
\left\{i+k : i \in \Lambda\right\}$.  A \embf{LIS} $f$ is
\embf{shift-invariant} or \embf{stationary} if
\[
f_{\Lambda+1}\left(\tau A \mid \tau \omega \right) \;=\;
f_{\Lambda}\left(A \mid \omega \right)
\]
for each $\Lambda \in \mathcal{S}_{b}$ and $\omega \in \Omega$.  We denote $\mathcal{G}_{\text{inv}}(f)$ the family of
shift-invariant chains consistent with a LIS $f$.

\begin{theo}[Ergodic chains]\label{th4}
Let $f$ be a shift-invariant LIS.
\begin{fleuvea}
\item A chain $\mu \in \mathcal{G}_{\text{inv}}(f)$ is extreme in
$\mathcal{G}_{\text{inv}}(f)$ if and only if $\mu$ is ergodic.
\item Let $\mu \in \mathcal{G}_{\text{inv}}(f)$.  If $\nu \in
\mathcal{P}_{\text{inv}}(\Omega,\mathcal{F})$ is such that $\nu \ll
\mu$, then $\nu \in \mathcal{G}_{\text{inv}}(f)$.
\item $\mathcal{G}_{\text{inv}}(f)$ is a face of
$\mathcal{P}_{\text{inv}}(\Omega,\mathcal{F})$. More precisely, if
$\mu, \, \nu \in \mathcal{P}_{\text{inv}}(\Omega,\mathcal{F})$ and
$0<s<1$ are such that $s\, \mu + (1-s) \, \nu \in
\mathcal{G}_{\text{inv}}(f)$ then $\mu, \, \nu \in
\mathcal{G}_{\text{inv}}(f)$.
\end{fleuvea}
\end{theo}

\section{Uniqueness results}
We shall prove two types of uniqueness results.  We start with the counterpart of a criterion proven by Georgii
(1974)\nocite{geo74} for measures determined by specifications.

\begin{theo}[One-sided boundary-uniformity]\label{th6}
Let $f$ be a LIS for which there exists a constant $c>0$ satisfying
the following property: For every $m\in\mathbb{Z}$ and every cylinder
set $A \in\mathcal{F}_{-\infty}^{m}$ there exists an integer $n<m$ such
that
\equa{\label{lis91}
f_{[n,m]}(A \mid \xi)\;\geq\; c \, f_{[n,m]}(A \mid \eta) \quad
\text{for all } \xi, \eta \in \Omega\;.
}
Then there exists at most one chain consistent with $f$.
\end{theo}

The main virtue of this criterion is its generality.  Existing uniqueness criteria (Harris, 1955\nocite{har55};
Iosifescu and Spataru, 1973\nocite{iosspa73}; Walters, 1955\nocite{wal75}; Berbee, 1987\nocite{ber87}; Stenflo,
2002\nocite{ste02}; Johansson and \"Oberg, 2002\nocite{johobe02}) require that the space $E$ have particular properties
(finite, countable, compact), and that the kernels satisfy appropriate non-nullness hypotheses.  Many of these criteria
are based on summability properties of the sequence of variations:
\eqna{\label{eq:rvar} \lefteqn{ \var{j}(f_{\{i\}}) \;\triangleq\;
}\nonumber\\ &&\sup\Bigl\{ \left| f_{\{i\}}(\xi_i\mid\xi_{-\infty}^i)
- f_{\{i\}}(\eta_i\mid\eta_{-\infty}^i) \right| : \xi, \eta \in
\Omega_{-\infty}^{i}, \, \xi_{j}^{i} = \eta_{j}^{i} \Bigr\}
}
for $j<i$.

\begin{prop}\label{cuc3} Assume that $E$ is a countable
set and $\mathcal{E}$ the discrete $\sigma$-algebra.  A LIS $f$
satisfies the one-sided boundary-uniformity condition \reff{lis91} if
it is \embf{uniformly non-null}:
\equa{\label{eq:unif}
\inf_{i\in\lat}\,\inf_{\omega\in\Omega_{\le i}}
f_{\{i\}}\bigl(\omega_i\bigm|\omega_{-\infty}^{i-1}\bigr) \;>\; 0\;,
}
and satisfies
\equa{\label{cuc5} \sup_{n\in\mathbb{Z}} \sum_{i \geq n}
\var{n} \left( f_{\{i\}}\right)\;<\;+\infty\;.
}
\end{prop}
We observe that when $f$ is stationary the last condition amounts to
summable variations: $\sum_{j < 0} \var{j} \left(
f_{\{0\}}\right)<+\infty$.

\medskip

Our second type of uniqueness result corresponds to the Dobrushin criterion for specifications.  The required
mathematical setting is richer.  We choose a bounded distance $d$ on $E$ and take $\mathcal{E}$ as the associated Borel
$\sigma$-algebra.  We endow $E^\lat$ with the product topology (so $\mathcal{F}$ is also Borel) and $\Omega\subset
E^{\lat}$ with the restricted topology.  The choice of distance is dictated by the type of measures to be analyzed.
For finite, or countable, alphabets the canonical choice is the discrete distance $d_{\text{disc}}(a,b) = 1$ if $a \neq
b$ and $0$ otherwise.

\begin{defi}
A LIS on $f$ on $(\Omega,\mathcal{F})$ is \embf{continuous} if the
functions $\Omega\ni\omega \;\longrightarrow\; f_{\Lambda} (A \mid
\omega )$ are continuous for all $\Lambda \in \mathcal{S}_{b}$ and all
$A\in\mathcal{F}_{\Lambda}$.
\end{defi}

In the case of specifications, continuity is associated with
Gibbsianness (non-nullness is also needed, see, e.g., the discussion
in Section 2.3.3 in van Enter, Fern\'andez and Sokal,
1993\nocite{vEFS_JSP}).  For $E$ finite, continuity is equivalent to
$\lim_{j\to-\infty}\var{j} \left( f_{\{i\}}\right)=0$.

\begin{rema}\label{rem:r1}
If the LIS $f$ is continuous and the space $\Omega$ is compact, then there always exists at least one compatible chain.
Indeed, the probability measures on a compact space form a (weakly) compact set. Hence, if
$(\Lambda_n)_{n\in\mathbb{N}}\subset\mathcal{S}_b$ is any exhausting sequence of regions and
$(\omega^{(n)})_{n\in\mathbb{N}}\subset\Omega$ any sequence of pasts, the sequence of measures
$f_{\Lambda_n}(\,\cdot\mid\omega^{(n)})$, $n\in\mathbb{N}$, has some accumulation point.  Continuity ensures that such
a limit belongs to $\mathcal{G}(f)$.  Therefore, for continuous LIS on a compact space of configurations, the following
theorems determine conditions for the existence of \emph{exactly} one compatible measure.
\end{rema}

For every $i \in \lat $ and every $\mathcal{F}_{\leq
i}$-measurable function $h$, the \embf{$d$-oscillation} of $h$ with
respect to the site $j \leq i$, is defined by
\equa{\label{dus9}
\delta_{j}^{d}(h) \;\triangleq\; \sup\left\{ \frac{\left| h(\xi) - h(\eta)
\right|}{d\left(\xi_{j},\eta_{j}\right)} : \xi, \eta \in
\Omega_{-\infty}^{i}\,, \, \xi \stackrel{\neq j}{=} \eta \right\}\;,
}
with the convention $0/0=0$ and where we introduced the notation
\equa{ \label{ccc5} \xi \stackrel{\neq j}{=} \eta
\quad \Longleftrightarrow \quad \xi_{i}=\eta_{i}\;, \; \forall \,
i\neq j}
(``$\xi$ equal to $\eta$ off $j$''). We introduce also the space of
functions of \emph{bounded $d$-oscillations}:
\equa{\label{eq:rlis.2}
\mathcal{B}_d \;\triangleq\; \Bigl\{
\,\mathcal{F}\mbox{-measurable } h:\sup_{j\in\lat}\delta_j^d(h)<\infty
\Bigr\} \;,
}
and its restrictions
\[
\mathcal{B}_d(\Lambda) \;\triangleq\; \Bigl\{ h\in\mathcal{B}_d :
h \,\mathcal{F}_\Lambda\mbox{-measurable } \Bigr\}
\]
for $\Lambda\subset\lat$.
%
The most general version of Dobrushin's strategy allows the use of a
``pavement'' of $\lat$ by finite intervals.  These intervals $V$ must be
chosen so that there is an appropriate control of the ``sensitivity''
of the averages $f_V$ to the configuration in $V_-$.

\begin{defi}[$d$-sensitivity estimator]\label{defest1}
Let $V \in \mathcal{S}_{b}$ and $f_V$ a probability kernel on
$\mathcal{F}_{\le m_V}\times \Omega$. A \embf{$d$-sensitivity
estimator} for $f_{V}$ is a nonnegative matrix $\alpha^{V} =\left(
\alpha_{ij}^{V}\right)_{i,j \in \lat}$ such that $\alpha^V_{ij}=0$ if
$i \notin V$ or $j \notin V_{-}$ and
\equa{\label{dus19}
\delta_{j}^{d} \left( f_{V}h \right)
\;\leq\; \sum_{i \in V}\delta_{i}^{d}(h) \, \alpha_{ij}^{V}
}
for all $ j \in V_{-}$ and $\mathcal{F}_V$-measurable functions $h
\in \mathcal{B}_d$.
\end{defi}

\begin{theo}[One-sided Dobrushin]\label{gencor1}
Let $f$ be a continuous LIS.  If there exist a countable partition
$\mathcal{P}$ of $\lat$ into finite intervals such that for each
$V \in \mathcal{P}$ there exists a $d$-sensitivity estimator
$\alpha^{V}$ for $f_{V}$ with
\equa{\label{gen30}
\sum_{j \in V_{-}} \alpha^V_{ij} \;<\; 1
}
for all $i\in\lat$, then there exists at most one chain consistent
with $f$.
\end{theo}

In particular, the partition can be trivial, namely $\mathcal{P}=\bigl\{\{i\}:i\in\lat\bigr\}$.  In the stationary
case, the estimators for such a partition are of the form $\alpha^{\{i\}}_{ij}=\alpha(i-j)$ for a certain function
$\alpha$ on the integers that takes value zero for non-positive integers. Dobrushin criterion becomes, then,
$\sum_{n\ge 1}\alpha(-n)<1$.
\smallskip

The customary way to construct $d$-sensitivity estimators for kernels
$f_V$ is resorting to the Vaserstein-Kantorovich-Rubinstein (VKR)
distance between measures on $\mathcal{F}_V$ for the distance
$d_{V}\left(\omega_{V},\sigma_{V}\right) \triangleq \sum_{i \in V}
d\left(\omega_{i},\sigma_{i}\right)$.  If we denote
$\stackrel{\circ}{f_{V}}$ the projection of each kernel $f_V$ over
$\Omega_V$:
\[
\stackrel{\circ}{f_{V}} \bigl( A \bigm|
\omega_{-\infty}^{l_V-1}\bigr) \;\triangleq\; f_{V}\bigl( \{
\sigma_{V} \in A \} \bigm| \omega_{-\infty}^{l_V-1} \bigr),
\quad \forall \; A \in \mathcal{F}_{V}, \;
\forall \; \omega_{-\infty}^{l_V-1} \in \Omega_{-\infty}^{l_V-1}\;,
\]
then the VKR distances between these projections are
\eqna{\label{vkr} \lefteqn{
\Bigl\|\stackrel{\circ}{f_{V}} \bigl( \,\cdot
\bigm| \xi_{-\infty}^{l_V-1}\bigr)\, - \stackrel{\circ}{f_{V}}
\bigl(\,\cdot \bigm| \eta_{-\infty}^{l_V-1}\bigr)\Bigr\|_{d_V}
\;=}\nonumber\\
&& \sup\left\{ \Bigl|\stackrel{\circ}{f_{V}} \bigl( h
\bigm| \xi_{-\infty}^{l_V-1}\bigr) - \stackrel{\circ}{f_{V}}
\bigl(h \bigm| \eta_{-\infty}^{l_V-1}\bigr)\Bigr| \;:\;
h\in\mathcal{B}_d(V)\,,\; {\rm osc}_V(h)\le 1 \right\}
\nonumber\\
\
}
where ${\rm osc}_V(h) = \sup\bigr\{|h(\sigma_V)-h(\omega_V)|/d_V(\sigma_V,\omega_V)\bigl\}$. Equivalently (see, for
instance, Dudley, 2002, Section 11.8\nocite{dud02}),
\eqna{\label{vkr2} \lefteqn{
\Bigl\|\stackrel{\circ}{f_{V}} \bigl( \,\cdot
\bigm| \xi_{-\infty}^{l_V-1}\bigr)\, - \stackrel{\circ}{f_{V}}
\bigl(\,\cdot \bigm| \eta_{-\infty}^{l_V-1}\bigr)\Bigr\|_{d_V}
\;=}\nonumber\\
&& \inf \left\{\int d(\sigma_V,\omega_V) \,  \rho(d\sigma_V,d\omega_V)
: \rho \in \mathcal{P}(\Omega \times \Omega)\right.\nonumber\\
&& \qquad\quad\left.\text{ with marginals }
\stackrel{\circ}{f_{V}} \bigl( \,\cdot \bigm|
\xi_{-\infty}^{l_V-1}\bigr) \text{ and }
\stackrel{\circ}{f_{V}} \bigr( \,\cdot \bigm|
\eta_{-\infty}^{l_V-1}\bigr)\right\}.
}
The VKR (canonical) $d$-estimator is defined by the coefficients
\equa{\label{int7}
C_{ij}^{V}(f) \;\triangleq\; \sup_{\substack{\xi, \eta \in
\Omega_{-\infty}^{l_{V} -1}\\ \xi \stackrel{\neq j}{=} \eta}}
\frac{\Bigl\| \stackrel{\circ}{f_{V}}\left( \, \cdot \mid
\xi \right) \,- \stackrel{\circ}{f_{V}}\left( \, \cdot \mid \eta
\right)\Bigr\|_{d_{V}}}{d\left(\xi_{j},\eta_{j}\right)}\;, \quad i
\in V, j\in V_-
}
and $C_{ij}^{V}(f)=0$ otherwise.

If the partition is trivial and $d$ is the discrete metric, each
$\Bigl\|\,\cdot\,\Bigr\|_{d_{\{i\}}}$ coincides with the variational
norm.  If the alphabet $E$ is countable, this means
\equa{\label{eq:del} C^{\{i\}}_{ij} \;=\; \delta_j(f_i) \triangleq\;\delta_j^{d_{{\rm disc}}}(f_i) \;=\; \sup\bigl\{
\left| f_i(\xi) - f_i(\eta) \right| : \xi, \eta \in \Omega \,, \, \xi \stackrel{\neq j}{=} \eta \bigr\}\;.}
and a sufficient condition for Dobrushin criterion \reff{lis91} is,
therefore,
\equa{\label{eq:rlis.1}
\sum_{j<i}\delta_j(f_i) \;<\;1\;,\quad
i\in\lat\;.
}
Besides the absence of non-nullness hypotheses, an advantage of Dobrushin criterion is that it determines a regime
where mixing properties can be determined, as we discuss in next section.

To conclude, we remark that in fact the two uniqueness criteria given in Theorems \ref{th6} and \ref{gencor1} give a
very strong form of uniqueness.
\begin{defi}[HUC]
  A LIS $f$ on $(\Omega,\mathcal{F})$ satisfies a \embf{hereditary
    uniqueness condition} (HUC) if for all intervals of the form
  $\Gamma=[k,+\infty[$, $k\in\mathbb{Z}$, and configurations
  $\omega\in\Omega$, the LIS $f^{(\Gamma,\omega)}$ defined by
\equa{\label{lis93}
f_{\Lambda}^{(\Gamma,\omega)}( \, \cdot \mid \xi)
\;=\; f_{\Lambda}(\, \cdot \mid \omega_{\Gamma_{-}} \,
\xi_{\Gamma})\;, \quad
\Lambda \in \mathcal{S}_{b}\,, \; \Lambda \subset \Gamma
}
admits at most one consistent unique chain.
\end{defi}

The two criteria given above involve bounds valid for all past conditions.  They remain, therefore, valid if only
particular pasts are considered as in \reff{lis93}.  This observation proofs the following corollary.

\begin{coro}\label{coro.r.1}
If a LIS satisfies the hypotheses of either Theorem \ref{th6} or Theorem \ref{gencor1}, then it also satisfies a HUC.
\end{coro}

We remark that, for similar reasons, the criteria of Harris (1955),
Stenflo (2002) and Johansson and \"Oberg (2002) also imply the
validity of a HUC.
\smallskip

\section{Results on loss of memory and mixing properties}
We place ourselves in the framework needed for the one-sided Dobrushin criterion ---$E$ with a topology defined by a
bounded metric $d$, $\mathcal{E}$ its Borel $\sigma$-algebra, $\Omega$ topologized with the restricted product
topology--- and take up all the related notions
---$d$-oscillations, functions of bounded oscillations, sensitivity
estimators.  To improve readability, we write the results only for a
trivial partition $\mathcal{P}$.  Versions for more general
partitions, of potential interest for coarse-graining
arguments, can be obtained in a straightforward manner from our
proofs by replacing sites by blocks of sites.

\begin{defi}
A \embf{$d$-sensitivity matrix} for a LIS $f$ is a matrix of the form
\equa{\label{eq:mat.r}
\alpha_{ij}\;\triangleq\; \begin{cases}
\alpha^{\{i\}}_{ij} &\mbox{ if } i>j\\
0 & \mbox{ otherwise}
\end{cases}
}
where each $\alpha^{\{i\}}_{ij}$ is a $d$-sensitivity estimator for
$f_i$, $i\in\lat$.
\end{defi}

\begin{theo}[Loss of memory]\label{memth1}
Let $f$ be a continuous LIS and $(\alpha_{ij})$ a $d$-sensitivity
matrix for $f$.  Then,
\begin{fleuvei}
\item For every $\Lambda\in\mathcal{S}_b$, $j < l_{\Lambda}$ and $h
  \in \mathcal{B}_d(\Lambda)$,
\equa{\label{eq:gen.r} \delta_{j}^{d}\left( f_{\Lambda}h \right) \;\leq\; \sum_{k \in \Lambda}
\delta_{k}^{d}(h) \left[ \sum_{l=1}^{|\Lambda|} \left( P_{\Lambda} \alpha \right)^{l}\right]_{kj} }
\item Assume that there exist a function $F : \lat^{2}
  \rightarrow \mathbb{R}^{+}$ satisfying the triangular inequality
  $F(i,j) \leq F(i,k) + F(k,j) \; \forall \, i,j,k \in \lat$
  such that
\equa{\label{mem5}
\gamma_{i} \;\triangleq\; \sum_{j < i} \alpha_{ij} \,
e^{F(i,j)}\;<\;1\;,
}
for each $i\in\lat$.  Then, for each $\Lambda \in \mathcal{S}_{b}$, $h \in \mathcal{B}_d(\Lambda)$ and $j <
l_{\Lambda}$.
\equa{\label{mem15} \delta_{j}^{d}\left( f_{\Lambda} h \right) \;\leq\; \frac{\gamma_{\Lambda}}{1-
\gamma_{\Lambda}} \sum_{k \in \Lambda} \delta_{k}^{d}(h) \, e^{-F(k,j)}\;, }
with $\gamma_{\Lambda} = \max_{i \in \Lambda} \gamma_{i}$.
\end{fleuvei}
\end{theo}

\begin{remas1}{2}
\item In the Markovian case $\alpha_{ij}=0$ if $|i-j|>1$.  Then
expression \reff{eq:gen.r} implies that for $h\in\mathcal{F}_{\{n\}}$
\equa{\label{eq:mar.r}
\delta^d_{-1}\Bigl(f_{[0,n]}(h)\Bigr) \;\le\; \gamma^n\,\delta^d_n(h)
}
with $\gamma=\sup_i\sum_j\alpha_{ij}$.  For $d$ discrete and estimators \reff{int7}, $\gamma$ is known as the
\embf{Dobrushin ergodic coefficient}.  If, in addition, $E$ is countable, $\Omega=E^\lat$ and $f$ shift-invariant, then
\equa{\label{eq.ddd.r}
\gamma\;=\;
1\,-\,\min_{\sigma_{-1},\omega_{-1}\in E} \sum_{\omega_0\in E}
f_{\{0\}}\bigl(\omega_0\bigm|\sigma_{-1}\bigr) \wedge
f_{\{0\}}\bigl(\omega_0\bigm|\omega_{-1}\bigr)\;.  }
\item If the alphabet $E$ is countable and the metric discrete we can
use the estimators \reff{eq:del}.  With this choice, \reff{mem15}
implies
\eqna{\label{eq:heu}
\lefteqn{ \delta_j(f_i) \;\le\, {\rm const}\,e^{-F(i,j)}
}\nonumber\\
&&\quad \Longrightarrow \quad
\delta_{-n}[f_{[0,m]}(A)] \;\le\; {\rm const}\, e^{-F(m,-n)}\;,\,
A\in \mathcal{F}_{\{m\}}\;.\quad
}
\end{remas1}

Published loss-of-memory results (Iosifescu, 1992\nocite{ios92};
Bressaud, Fern\'andez and Galves, 1999\nocite{bfg2}) resort instead to
the variations \reff{eq:rvar}.  Comparisons can only be made through
the obvious inequalities
\[
\delta_j[f_{\{i\}}(h)]\;\le\;\var{j}[f_{\{i\}}(h)]
\;\le\;\sum_{k\le j}\delta_k[f_{\{i\}}(h)]\;.
\]
For LIS with an exponentially decaying dependence on the past,
\reff{eq:heu} implies an exponential loss of memory with an identical
rate, in terms either of oscillations or of variations.  This should
be contrasted with the results in Bressaud, Fern\'andez and Galves
(1999)\nocite{bfg2} where there is an infinitesimal loss of rate.  LIS
with a power-law dependence can be treated by taking
$F(i,j)=c\log(1+|i-j|)$.  In terms of variations, the loss of memory
implied by \reff{eq:heu} is also a power law but with a power
decreased by one unit.  Bressaud, Fern\'andez and Galves
(1999)\nocite{bfg2} obtain, instead, the same power.

Furthermore, it is relatively simple to construct examples falling
outside the scope of all preexisting loss-of-memory results, but for
which Theorem \ref{memth1} applies.  Consider, for instance, the
2-letter alphabet $E=\{0,1\}$ and a shift-invartiant LIS defined by
singletons
\equa{\label{eq:exa.r}
f\bigl(\omega_0=1\bigm| \omega_{-\infty}^{-1}\bigr) \;=\;
\sum_{i\le 0} a_i\,\omega_i\;,
}
for a sequence $\{a_i\}_{i\le 0}$ of non-negative numbers.  The
estimators \reff{eq:del} yield a sensitivity matrix
\equa{\label{eq:exa2.r} \alpha_{ij} \;=\; \delta_j\left(f_{\{i\}}\right) \;=\; a_{i-j} }
for $i>j$, and zero otherwise.  Theorem \reff{memth1} is therefore
applicable as long as $\sum_{i\le 0}a_i<1$.  On the other hand, for each
$0<\varepsilon<1$, the choice
\equa{\label{eq:exa3.r}
a_{-k} \;=\;
\frac{1-\varepsilon}{M_\varepsilon}\,\frac{1}{k^{1+\varepsilon}}
}
with $M_\varepsilon=\sum_{k\ge 1} k^{-(1+\varepsilon)}$, satisfies
\[
\var{j}(f_{\{i\}}) \;\ge\; \frac{1}{(i-j-1)^\varepsilon}
\]
for $i-j \ge 2$.  Thus, this LIS is not covered by the results of Iosifescu (1992) or of Bressaud, Fern\'andez and
Galves (1999).  It also does not satisfy any uniqueness criteria except one-sided Dobrushin's.
\medskip

The following mixing results form the LIS version of a well known
chapter in the theory for Gibbs measures (see, for example, chapter V
in Simon, 1993\nocite{sim93}).  Their proofs, presented in Section
\ref{s.plmm}, follow the guidelines of the statistical mechanical
proofs.  They require a compact $\Omega$.  We observe that example
\reff{eq:exa.r}--\reff{eq:exa3.r} shows that our results are
complementary to those existing in the literature, which are based on
variations rather than oscillations (Bressaud, Fern\'andez and Galves,
1999, and references therein).

\begin{theo}\label{comth1}
  Assume $\Omega$ compact and let $f$ and $\widetilde{f}$ be two LIS
  on $(\Omega, \mathcal{F})$ with $f$ continuous and with a unique
  consistent measure. Assume also that for each $i \in \lat $ there
  exists a measurable function $b_{i}$ on $\Omega$ such that
\equa{\label{com5}
\Bigl\| \stackrel{\circ}{f_{\{i\}}}(\, \cdot \mid \omega)\, -
\stackrel{\circ}{\widetilde{f}_{i}}(\, \cdot \mid \omega) \Bigr\|_{d}
\;\leq\; b_{i}(\omega)
}
for every configuration $\omega \in \Omega_{-\infty}^{i-1}$.  Then,
for all $\mu \in \mathcal{G}(f)$, $\widetilde{\mu} \in \mathcal{G}
(\widetilde{f})$ and $\Lambda \in\mathcal{S}_b$
\equa{\label{com7}
\bigl| \mu(h) - \widetilde{\mu}(h) \bigr| \;\leq\;
\sum_{k \in \Lambda \cup \Lambda_{-}} \widetilde{\mu} \left( b_{k}
\right)\, \delta_{k}^{d}\left( f_{[k+1, m_{\Lambda}]}h\right)
}
for every $h \in \mathcal{B}_d(\Lambda)$.
\end{theo}

Let us denote $D \triangleq \sup_{x,y \in E} d(x,y)$ and for a measure
$\mu$ on $\mathcal{F}$ and $\mathcal{F}$-measurable functions $h_{1}$
and $h_{2}$
\[
\text{Cor}_{\mu}\left(h_{1},h_{2}\right) \;\triangleq\; \Bigl| \mu
\left(h_{1} \, h_{2} \right) - \mu(h_{1}) \mu(h_{2})\Bigr|\;.
\]
\begin{theo}\label{corth1}
  Assume $\Omega$ compact and let $f$ be a LIS on $(\Omega,
  \mathcal{F})$ that is continuous and and with a unique consistent
  measure.  Let $\mu$ be the unique probability measure in
  $\mathcal{G}(f)$.  Then for every $\Lambda, \Delta \in
  \mathcal{S}_{b}$ such that $m_{\Delta} < l_{\Lambda}$,
\equa{\label{cor5}
\text{Cor}_{\mu}\left(h_{1},h_{2}\right) \;\leq\;
\frac{D^{2}}{4}\sum_{k \leq m_{\Delta}} \delta_{k}^{d}
\left( f_{[k+1,m_{\Lambda}]}h_{1} \right) \delta_{k}^{d} \left(
f_{[k+1,m_{\Delta}]} h_{2}\right)
}
for all functions $h_{1} \in \mathcal{B}_d(\Lambda)$ and $h_{2} \in
\mathcal{B}_d(]-\infty, m_{\Delta}]).$
\end{theo}

Next corollary offers a more quantitative consequence of this theorem.
For all $\Lambda \in \mathcal{S}_{b}$ we define the
$\Lambda$-projection
\[
\left( P_{\Lambda} \right)_{kj} \;=\; \begin{cases} 1 & \text{ if }
k=j \text{ and } k \in \Lambda\\
0 & \text{ otherwise\;.} \end{cases}
\]
For a matrix $\left( A_{kj}\right)_{k,j \in \lat }$ with
nonnegative entries, we denote
\equa{\label{gen7}
\left[ \frac{A}{1 - A} \right]_{kj} \;\triangleq\; \sum_{n \geq 1}
\left[A^{n} \right]_{kj}\;.
}
These are well-defined sums on $[0,+\infty]$.

\begin{coro}\label{corcor1}
Consider the hypotheses of the previous theorem and let
$(\alpha_{ij})$ be a $d$-sensitivity matrix for $f$.
\begin{fleuvei}
\item If $ h_{1} \in \mathcal{B}_d(\Lambda)$ and
$ h_{2} \in \mathcal{B}_d(]-\infty, m_{\Delta}])$,
\equa{\label{cor29}
\text{Cor}_{\mu}\left(h_{1},h_{2}\right) \;\leq\;
\frac{D^{2}}{4}\sum_{k \leq m_{\Delta}}\sum_{l \in \Lambda}
\delta_{l}^{d} (h_{1}) \left[ \frac{P_{\Lambda}
\alpha} {1-P_{\Lambda} \alpha} \right]_{lk} \delta_{k}^{d} \left(
f_{[k+1,m_{\Delta}]} h_{2}\right)\;.
}
\item If $ h_{1} \in \mathcal{B}_d(\Lambda)$ and
$h_{2} \in \mathcal{B}_d(\Delta)$,
\equa{\label{cor31}
\text{Cor}_{\mu}\left(h_{1},h_{2}\right) \;\leq\;
\frac{D^{2}}{4} \sum_{l \in \Delta} \sum_{m \in \Lambda}
\delta_{m}^{d} (h_{1})\, \delta_{l}^{d} (h_{2})\, A_{ml}\;,
}
where
\[
A_{ml} \;\triangleq\; \left[ \frac{P_{\Lambda}
\alpha}{1-P_{\Lambda} \alpha} \right]_{ml} +
\sum_{k \leq m_{\Delta}} \left[
\frac{P_{\Lambda} \alpha}{1- P_{\Lambda} \alpha} \right]_{mk}
\left[ \frac{P_{[k+1,m_\Delta]}\, \alpha}{1 - P_{[k+1,m_\Delta]}\, \alpha}
\right]_{lk} \;.
\]
\end{fleuvei}
\end{coro}

The following proposition is useful to estimate the different matrices
appearing in this corollary.

\begin{prop}\label{p.ff1}
If $(\alpha_{ij})$ is a matrix satisfying \reff{mem5}, then for each
$\Lambda\in\mathcal{S}_b$
\equa{\label{mem21}
\left[ \frac{P_{\Lambda}\,\alpha}{1-P_{\Lambda}\, \alpha}\right]_{kj}
\;\leq \;
\frac{\gamma_{\Lambda}}{1-\gamma_{\Lambda}} \; e^{-F(k,j)}\;.
}
\end{prop}

\section{Proofs for the general framework}\label{s.pgf}

\subsection{Singleton consistency for chains}
The fact that the objects defined by \reff{liss9} are kernels from
$\mathcal{F}_{\leq m_\Lambda} \times \Omega$ to the interval $[0,1]$
follows immediately from the properties of the kernels $f_i$.  Their
normalization is proven by induction, using the fact that
\[
f_{\{i\}}(1\mid\cdot\,)\;=\;f_{\{i\}}\bigl(\Omega_{\le i}\bigm| \cdot\,\bigr)\;=\;1
\]
and the inductive step
\[
f_\Lambda\bigl(\Omega_{\le m_\Lambda}\bigm|\omega\bigr) \;=\;
f_{[l_\Lambda,m_\Lambda-1]}\Bigl( \bigl(f_{m_\Lambda}\bigl(\Omega_{\le
m_\Lambda}\bigr) \Bigm| \omega\Bigr) \;=\;
f_{[l_\Lambda,m_\Lambda-1]}\bigl( 1 \bigm|\omega\bigr)\;=\; 1\;,
\]
for $\omega\in\Omega_{\le l_\Lambda}$.
\smallskip

Properties (a) and (b) of the definition \ref{lis1} of LIS are an
immediate consequence of similar properties of the kernels $f_i$.  To
prove consistency, we first remark that for $l\le m\le p$,
$\omega\in\Omega$ and any $\mathcal{F}_{\le p}$-measurable function
$h$,
\eqna{\label{eq:uf1}
\bigl(f_{[l,m]}\, f_{[l,p]}\bigr)(h\mid\omega) &=&
f_{[l,m]} \Bigl( f_{[l,p]}(h)\Bigm| \omega\Bigr) \nonumber\\
&=&  f_{[l,p]}(h\mid \omega)\;  f_{[l,m]}(1\mid \omega)\nonumber\\
&=& f_{[l,p]}(h\mid \omega)\;.
}
The second equality is due to the proven property (b) of Definition
\ref{lis1} plus the fact that $f_{[l,p]}\left( h \mid \cdot \,\right)$
is $\mathcal{F}_{\leq l-1}$-measurable.  The last equality is the just
proven normalization.  Identity \reff{eq:uf1} justifies the last
equality in the following string of identities, valid for $l\le m<
p$,
\equa{\label{eq:uf2}
f_{[l,p]}\, f_{[l,m]} \;=\; f_{[l,m]} \,f_{[m+1,p]} \, f_{[l,m]}
\;=\; f_{[l,m]}\, f_{[l,p]} \;=\; f_{[l,p]} \;.
}
The other equalities are simply due to definition \ref{liss9}.  A
similar identity is trivially true for $l\le m=p$.  Consistency
follows for, if $\Delta \supset \Lambda$:
\[
f_{\Delta} \, f_{\Lambda}\;=\; f_{[l_\Delta,l_\Lambda-1]}\,
f_{[l_\Lambda,m_\Delta]}\, f_{[l_\Lambda,m_\Lambda]} \;=\;
f_{[l_\Delta,l_\Lambda-1]}\, f_{[l_\Lambda,m_\Delta]}
\; =\; f_{\Delta}\;.
\]
We used \reff{eq:uf2} in the middle identity and we assumed
$l_\Delta<l_\Lambda$, otherwise we revert to \reff{eq:uf2}.

The remainder of the proof relies on the following observation valid
for any measure $\mu$ on $\mathcal{F}$ and
any$\Lambda\in\mathcal{S}_b$:
\equa{\label{eq:uf3}
\mu \,f_i=\mu\,,\; \forall i\in\Lambda \quad\Longrightarrow\quad
\mu\,f_\Lambda \;=\; \mu\;.
}
This is proven by induction on the cardinality of $\Lambda$ through the identity
\[
\mu \, f_{\Lambda} \;=\; \mu \, f_{l_\Lambda} \,
f_{[l_\Lambda+1,m_\Lambda]} \;=\;  \mu \,
f_{[l_\Lambda+1,m_\Lambda]}\;.
\]

Property \reff{eq:uf3} directly proves the non-trivial inclusion in
\reff{liss7}. Furthermore, it yields uniqueness.  Indeed, consider a
LIS $\left(g_{\Lambda}\right)_{\Lambda \in \mathcal{S}_{b}}$
consistent with the family $\left(f_{i}\right)_{i \in\lat}$.  By
\reff{eq:uf3} $g_\Lambda$ must be consistent with $f_{\Lambda}$ for
each $\Lambda\in\mathcal{S}_b$.  But then, if $\omega\in\Omega$ and
$h$ is $\mathcal{F}_{\le m_\Lambda}$-measurable
\[
g_{\Lambda}\bigl(h\bigm|\omega\bigr)\;=\;
g_{\Lambda}\Bigl(f_{\Lambda}(h)\Bigm|\omega\Bigr) \;=\;
f_{\Lambda}\bigl(h\bigm|\omega\bigr)\,g_{\Lambda}\bigl(1\bigm|\omega\bigr)
\;=\;f_{\Lambda}\bigl(h\bigm|\omega\bigr)\;.
\]
The second identity is a consequence of the
$\mathcal{F}_{l_\Lambda-1}$-measurability of $f_\Lambda(h|\cdot\,)$
plus property (b) of Definition \ref{lis1}.  The last equality is the
normalization of $g_\Lambda$.$\quad \Box$

\subsection{Extreme chains}
We start with general results on probability kernels.
\begin{prop}\label{prop1}
Let $\mathcal{B}$ be a sub-$\sigma$-algebra of $\mathcal{F}$, $\pi$ a
probability kernel on $\mathcal{B} \times \Omega$ and $\mu \in
\mathcal{P}(\Omega, \mathcal{F})$ such that $\mu \pi = \mu$ on
$\mathcal{B}$. Then:
\begin{fleuvei}
\item  The system
\equa{\label{eq:ex1.r}
\mathcal{I}_{\pi}^{\mathcal{B}}(\mu) \;\triangleq\;
\Bigl\{ A \in \mathcal{B}: \pi(A \mid \cdot \,) = \ind{A}(\, \cdot
\,) \; \mu\text{-a.s.} \Bigr\}
}
is a $\sigma$-algebra.
\item For all $\mathcal{B}$- measurable functions $h : \Omega
\rightarrow [0, + \infty[\,$,
\equa{\label{eq:ex2.r}
(h \mu)\,\pi= h \mu \text{ on } \mathcal{B}\quad \text{if and only if}
  \quad h \text{ is } \mathcal{I}_{\pi}^{
\mathcal{B}}(\mu) \text{-measurable\;.}
}
\end{fleuvei}
\end{prop}
\textbf{Proof}\\
{\bf (i)} Clearly $\Omega \in \mathcal{I}_{\pi}^{ \mathcal{B}}(\mu)$. For each $A \in
\mathcal{I}_{\pi}^{\mathcal{B}}(\mu)$,
\[
\pi(A^{c} \mid \cdot \, ) \;=\; 1- \pi(A \mid \cdot \, )
\;=\; 1 - \ind{A} \ (\mu\text{-a.s.}) \;=\; \ind{A^{c}} \
(\mu\text{-a.s.})\;.
\]
Likewise, for each sequence $(A_{n})_{n \in \mathbb{N}}$ of disjoint
sets in $\mathcal{I}_{\pi}^{\mathcal{B}}(\mu)$,
\[
\pi\bigl(\cup A_{n} \bigm| \cdot \, \bigr) \;=\;
\sum_{n \in \mathbb{N}}\pi(A_{n} \mid \cdot \, )\;=\;
\sum_{n \in \mathbb{N}} \ind{A_{n}} \ (\mu\text{-a.s.})
\;=\; \ind{\cup A_{n}} \ (\mu\text{-a.s.})\;.
\]
Finally, if $A,B \in
\mathcal{I}_{\pi}^{\mathcal{B}}(\mu)$, then
\[
\pi(A \cap B \mid \cdot \, ) \;\leq\;
\pi(A \mid \cdot \, ) \wedge \pi(B \mid \cdot \, )
\;=\; \ind{A} \wedge \ind{B} \ (\mu\text{-a.s.})
\;=\; \ind{A \cap B} \ (\mu \text{-a.s.})
\]
and, by the consistency of $\mu$ with $\pi$,
\[
\mu\Bigl(\ind{A \cap B}- \pi(A \cap B \mid \cdot \, )\Bigr)
\;=\; \mu(A \cap B)- \mu \pi(A \cap B)\;=\;0\;.
\]
Thus
\[ \pi(A \cap B)= \ind{A \cap B} \; \mu\text{-a.s.}\]
\smallskip

\noindent {\bf (ii)}
Let us assume that $(h \mu)\,\pi=h \mu$ on $\mathcal{B}$. To prove
necessity it suffices to show that $\{h \geq c\}
\in \mathcal{I}_{\pi}^{ \mathcal{B}}(\mu),$ for all $c>0.$ Let us
fix some $c>0$ and denote $g=\ind{h \geq c}$. We have
\begin{align*}
\mu\Bigl((1-g)\, h \, \pi(g) \Bigr) \;& =\;
 (h \, \mu)\bigl( \pi (g)\bigr) - \mu \bigl(g \, h \, \pi (g)\bigr)
\;=\; (h \, \mu)(g) - \mu \bigl(g \, h \, \pi (g)\bigr)\\
& =\; \mu\Bigl(g \, h \, \bigl(1- \pi (g)\bigr)\Bigr)\;.
\end{align*}
But $g h \geq c g$ and $1- \pi(g) \geq 0$, hence
\begin{align*}
\mu\Bigl((1-g)\, h \, \pi(g) \Bigr)\; & \geq\; c \,
\mu\Bigl(g\bigl(1 - \pi(g)\bigr)\Bigr) \;=\; c\, \mu\bigl(\pi(g)\bigr)
- c \, \mu\bigl(g \, \pi(g)\bigr)\\
& =\; c\, \mu\Bigl( (1-g) \, \pi(g)\Bigr)\;.
\end{align*}
We obtain that $\mu \left( \ind{\{h<c\}} \, (h-c) \, \pi(g) \right) \geq 0$, which implies $ \ind{\{h<c\}} \pi(g) = 0
\; \; \mu$-a.s. Therefore,
\[
\pi(g)\;=\; g \, \pi(g) + \ind{\{h<c\}}\, \pi(g) \;\leq\; g
\;\; \mu \text{-a.s.}
\]
Furthermore, $\mu\bigl(g- \pi(g)\bigr)= 0$ by the consistency of $\mu$
with $\pi$.  This fact, together with the previous inequality, allow
us to conclude that $\pi(g) =g$ $\mu$-a.s., that is $\{h \geq c\} \in
\mathcal{I}_{\pi}^{\mathcal{B}}(\mu)$.
\smallskip

Conversely, assume that $h$ is
$\mathcal{I}_{\pi}^{\mathcal{B}}(\mu)$-measurable. By the standard
machinery of measure theory sufficiency follows if we show
for all $A \in \mathcal{I}_{\pi}^{\mathcal{B}}(\mu)$ that
$(\ind{A}\,\mu)\pi = \ind{A} \mu \text{ on } \mathcal{B}$. If $ B \in
\mathcal{B}$,
\begin{align*}
(\ind{A}\, \mu)\,\pi(B)\; & =\; (\ind{A}\, \mu )\,\pi(A
\cap B) + (\ind{A}\, \mu)\,\pi(B \setminus A)\\
& \leq\; \mu \pi(A \cap B) + (\ind{A}\, \mu)\,\pi\left(A^{c}\right)\;.
\end{align*}
The consistency of $\mu$ with $\pi$ implies that the second term of
the last line is zero.  Thus we have proved that
\equa{\label{ext3}
(\ind{A} \mu)\,\pi(B) \;\le\; (\ind{A} \mu)(B)\;.
}
By the same token,
\equa{\label{ext7}
(\ind{A} \mu)\,\pi \left(B^{c}\right) \;\leq\; (\ind{A} \mu
)\left(B^{c}\right)\;.
}
But the consistency of $\mu$ with $\pi$ implies that the sum of
the LHS of \reff{ext3} and \reff{ext7} equals the sum of the
corresponding RHS, namely $\mu(A)$.  We conclude that
$ (\ind{A} \mu)\,\pi(B) = (\ind{A} \mu)(B)$.$\quad \Box$
\begin{coro}\label{corext1}
Let $\Pi$ be a non-empty set of probability kernels $\pi$ defined on
$\mathcal{F}_{\pi} \times \Omega$, where $\mathcal{F}_{\pi}$ is a
sub-$\sigma$-algebra of $\mathcal{F}$. Let us denote
\equa{\label{eq:pgf3}
\mathcal{G}(\Pi) \;=\; \Bigl\{\mu \in \mathcal{P}(\Omega,\mathcal{F})
: \mu \, \pi = \mu \text{ on } \mathcal{F}_{\pi}
\text{ for all } \pi \in \Pi\Bigr\}
}
and for each $\mu \in \mathcal{G}(\Pi)$,
\equa{\label{eq:pgf4}
\mathcal{I}_{\Pi}(\mu) \;=\; \bigcap_{\pi\in \Pi}
\mathcal{I}_{\pi}^{\mathcal{F}_{\pi}}(\mu)
}
be the $\sigma$-algebra of all $\mu$-almost surely $\Pi$-invariant
sets. Then $\mu$ is trivial on $\mathcal{I}_{\Pi}(\mu)$ if $\mu$ is
extreme in $\mathcal{G}(\Pi)$.
\end{coro}
\textbf{Proof}\\
Suppose $\mu$ is not trivial on $\mathcal{I}_{\Pi}(\mu)$ and take
$A \in \mathcal{I}_{\Pi}(\mu)$ such that $0< \mu(A) <1$.  The measures
\[
\nu\;=\; \mu( \, \cdot \mid A) \;\triangleq\; h\, \mu \quad \text{with
} h =\frac{\ind{A}}{\mu(A)}
\]
and
\[
\nu'\;=\; \mu( \, \cdot \mid A^{c}) \;\triangleq\; h' \mu \quad
\text{with } h' =\frac{\ind{A^{c}}}{\mu(A^{c})}
\]
satisfy $ \nu \neq \nu'$ and $\mu = \mu(A) \, \nu + \mu(A^{c}) \, \nu'$.  The functions $h$ and $h'$ are
$\mathcal{I}_{\pi}^{\mathcal{F}_{\pi}}(\mu)$-measurable, for all $\pi \in \Pi$. Thus, (ii) of Proposition \ref{prop1}
implies that $\nu, \, \nu' \in \mathcal{G}(\Pi)$, a fact that contradicts the extremality of $\mu. \quad \Box$
\begin{lemm}\label{lemext1}
Let $f$ be a LIS defined on $(\Omega, \mathcal{F})$ and $\mu \in
\mathcal{G}(f)$. Let us denote $\mathcal{F}_{-\infty}^{\mu}$ the
$\mu$-completion of $\mathcal{F}_{-\infty}$. Then
\equa{\label{lemext2}
\bigcap_{n \geq 0}\mathcal{I}_{f_{[k-n,k]}}^{\mathcal{F}_{\leq
    k}}(\mu) \;=\; \mathcal{F}_{-\infty}^{\mu}
}
for each $k \in \lat$ and
\equa{\label{lemext3}
\bigcap_{\Lambda \in \mathcal{S}_{b}}
\mathcal{I}_{f_{\Lambda}}^{\mathcal{F}_{\leq
m_\Lambda}}(\mu) \;=\; \mathcal{F}_{-\infty}^{\mu}\;.
}
\end{lemm}
\textbf{Proof}\\
Identity \reff{lemext2} follows from the observation that for each $B
\in \bigcap_{n}\mathcal{I}_{f_{[k-n,k]}}^{\mathcal{F}_{ \leq k}}(\mu)$
the set $A \triangleq \bigcap_{n} \left\{ f_{[k-n,k]}(B \mid \cdot
\,)=1 \right\}$ satisfies $A=B \; \mu\text{-a.s.}$ and $A \in
\mathcal{F}_{-\infty}$. Equality \reff{lemext3} is a consequence of
\reff{lemext2} because
\[
\displaystyle{\bigcap_{\Lambda \in
\mathcal{S}_{b}} \mathcal{I}_{f_{\Lambda}}^{\mathcal{F}_{\leq
m_\Lambda}}(\mu) \;=\; \bigcap_{k \in \mathbb{Z}} \bigcap_{n \geq 0}
\mathcal{I}_{f_{[k-n,k]}}^{\mathcal{F}_{\leq k}}(\mu)}. \quad \Box
\]
\medskip

\noindent\textbf{Proof of Theorem \protect \ref{th2}}\\
\noindent{\bf (a)} It is immediate.
\smallskip

\noindent{\bf (b) $(\Rightarrow)$} The implication follows readily
from Corollary \ref{corext1} and the fact that, by \reff{lemext3},
$\bigcap_{\Lambda \in \mathcal{S}_{b}}
\mathcal{I}_{f_{\Lambda}}^{\mathcal{F}_{\leq m_\Lambda}}(\mu)$ is
$\mu$-trivial if and only if $\mu$ is trivial on
$\mathcal{F}_{-\infty}$.
\smallskip

\noindent{\bf (c) $(\Rightarrow)$}  Let $\mu, \; \nu \in \mathcal{G}(f)$ such that $\nu \ll \mu$.  There exists a
$\mathcal{F}$-measurable non-negative function $g$ such that
\[\nu \;=\; g\, \mu\;.\]
Let us consider, for each $k \in \lat \; \mu_{k} \triangleq \mu
\big\vert_{\mathcal{F}_{\leq k}}$ and $\nu_{k} \triangleq \nu
\big\vert_{\mathcal{F }_{\leq k}}$. As in particular $\nu_{k} \ll
\mu_{k}$ on $\mathcal{F}_{\leq k}$, there
exists $g_{k} \geq 0$, $\mathcal{F}_{\leq k}$-measurable, satisfying
$ \nu_{k}= g_{k} \, \mu_{k}$ on $\mathcal{F}_{\leq k})$.  All we have
to prove is that
\equa{\label{ext19}
g_{k} \text{ is } \mathcal{F}_{-\infty}^{\mu}\text{-measurable}
\quad \forall \, k \in \lat\;.
}
Indeed, by the reverse martingale theorem $g_{k} = g$
$\mu$-a.s. Therefore, $g$ inherits the
$\mathcal{F}_{-\infty}^{\mu}$-measurability and, thus, it is
$\mu$-a.s.\ equal to a $\mathcal{F}_{-\infty}$-measurable function.

To prove \reff{ext19} we observe that
since $\nu \in \mathcal{G}(f),$
\[
g_{k}\, \mu_{k}\, f_{[k-n,k]}\;=\; g_{k} \,\mu_{k}
\]
on $\mathcal{F}_{ \leq k}$ for all $n \in \mathbb{N}$. As $g_{k}$ is
$\mathcal{F}_{\leq k}$-measurable, we conclude from Proposition
\ref{prop1} that
$ g_{k} \text{ is }
\bigcap_{n}\mathcal{I}_{f_{[k-n,k]}}^{\mathcal{F}_{\leq k}}(\mu)
\text{-measurable}$.  Its, $\mathcal{F}_{-\infty}^{\mu}$-measurability
follows, hence, from \reff{lemext2}.
\smallskip

\noindent{\bf (b) $(\Leftarrow)$} Assume $\mu$ is a trivial measure on $\mathcal{F}_{-\infty}$ and suppose that there
exist $s: 0<s<1$ and $\nu, \nu' \in \mathcal{G}(f)$ such that $\mu=s \, \nu + (1-s) \, \nu'$. As $\nu,\nu' \ll \mu$, by
(c) ($\Rightarrow$) there exist $\mathcal{F}_{-\infty}$-measurable functions $h, h'\geq 0$ such that $\nu = h \mu$ and
$\nu' = h' \mu$.  But the triviality of $\mu$ on $\mathcal{F}_{-\infty}$ implies that $ h=h'=1$ $\mu$-a.s. Thus $ \mu =
\nu = \nu'$.
\smallskip

\noindent{\bf (c) $(\Leftarrow)$} This is an immediate consequence of
Proposition \ref{prop1} plus the fact that $h$ is
$\mathcal{I}_{f_{\Lambda}}^{\mathcal{F} \leq
  m_\Lambda}(\mu)$-measurable for all
$\Lambda \in \mathcal{S}_{b}$.
\smallskip

\noindent{\bf (d)} Let $\mu, \nu \in \mathcal{G}(f)$ such that $\mu=
\nu$ on $\mathcal{F}_{-\infty}.$ Consider $\widetilde{\mu} \triangleq
\frac{1}{2} \, \mu + \frac{1}{2} \, \nu \; \in \mathcal{G}(f)$. Since
$\mu \ll \widetilde{\mu}$ and $\nu \ll \widetilde{\mu},$ assertion (b)
implies that
$ \mu = f \widetilde{\mu}$ and $\nu = g \widetilde{\mu}$ for
$\mathcal{F}_{-\infty}$-measurable functions $f$ and $g$.
But $\mu=\nu=\widetilde{\mu}$ on $\mathcal{F}_{-\infty}$, so $f=g$
$\mu$-a.s. and therefore $\mu = \nu.$
\smallskip

\noindent{\bf (e)} It is an immediate consequence of (b) and
(d). $\quad \Box$

\subsection{Triviality and short-range correlations}
The proofs involve standard arguments.  We include them for
completeness.
\smallskip

\noindent\textbf{Proof of Proposition \protect\ref{theo2}}\\
\noindent{\bf (a) $\Rightarrow$ (c)} Let $A \in \mathcal{F}$ and
$k \in \lat$. Since $ \mathcal{F}_{-\infty} = \bigcap_{n \geq
1} \mathcal{F}_{\leq k-n}$, the reverse martingale theorem yields
\equa{\label{triv5}
\mu \bigl( A \bigm| \mathcal{F}_{\leq k-n} \bigr)
\;\tend{n}{+\infty}{L^{1}(\mu)}\;
\mu \bigl( A \bigm| \mathcal{F}_{-\infty} \bigr)\;.
}
The assumed triviality of $\mu$ on $\mathcal{F}_{-\infty}$ implies
that $\mu \left( A \mid \mathcal{F}_{-\infty}\right)=
\mu(A) \quad \mu\text{-a.s.}$  We deduce that for each $\varepsilon
>0$, there exists $\Delta \in \mathcal{S}_{b}$ such that
\equa{\label{triv13} \mu \Bigl( \bigl| \mu \left(A \mid
  \mathcal{F}_{\Delta_{-}} \right) - \mu(A) \bigr|\Bigr)
\;< \; \varepsilon\;.
}
Hence, for all $\Lambda \in \mathcal{S}_{b} : \Lambda \supset \Delta,$
\begin{eqnarray*}
\sup_{B \in \mathcal{F}_{\Lambda_{-}}}
\Bigl| \mu(A \cap B) - \mu(A)\,\mu(B)\Bigr|  & \leq & \sup_{B \in
\mathcal{F}_{\Delta_{-}}} \Bigl| \mu(A \cap B) - \mu(A)\,\mu(B)\Bigr|\\
&=& \Bigl|\mu \Bigl( \bigr[ \mu \left( A \mid
\mathcal{F}_{\Delta_{-}}\right) - \mu(A)  \bigr]\,\ind{B}
\Bigr)\Bigr| \\
& \leq& \mu \Bigl( \bigl| \mu \left( A \mid
\mathcal{F}_{\Delta_{-}}\right) - \mu(A)  \bigr| \Bigr)\\
& < & \varepsilon.
\end{eqnarray*}

\noindent{\bf (b) $\Rightarrow$ (a)} Fix $B \in \mathcal{F}_{-\infty}$
and consider $\mathcal{D} \triangleq \left\{ A \in \mathcal{F} : \mu(A
\cap B) = \mu(A)\mu(B) \right\}$. It is straightforward to see that
$\mathcal{D}$ is a $\lambda$-system. By assumption $\mathcal{D}$
contains all cylinder events, so $\mathcal{D}= \mathcal{F}$ [Dynkin's
$\pi$-$\lambda$ theorem]. In particular $B \in \mathcal{D}$, thus
$\mu(B) = \left( \mu(B)\right)^{2}$ and thereby $\mu(B) = 0$ or
$1$. $\quad \Box$
\medskip

\noindent\textbf{Proof of Theorem \protect\ref{th3}}\\
\noindent{\bf (a)} Let $h$ be a bounded local function on $\Omega$. As
$\mu$ is consistent with $f$, $f_{\Lambda_{n}}h$ coincides with
$\mu\left(h \mid \mathcal{F}_{(\Lambda_{n})_{-}}\right)$, $\mu$-a.s.,
for $n$ sufficiently large.  Therefore, by the reverse martingale
convergence theorem we conclude that
\[
f_{\Lambda_{n}}h \;\tend{n}{+\infty}{\ }\;
\mu\left( h \mid \mathcal{F}_{-\infty}\right) \quad \mu\text{-a.s.}
\]
This implies assertion (a) because $\mu$ is trivial on
$\mathcal{F}_{-\infty}$.
\smallskip

\noindent{\bf (b)} It is a consequence of assertion (a) and the fact
that if $\Omega$ is compact and metric, the space of local continuous
functions on $\Omega$ contains a countable subset which is dense with
respect to the uniform-norm.  $\quad \Box$

\subsection{Ergodicity}
We need a well known result of ergodic theory.  See, for instance,
Georgii (1988), Theorem 14.5, for a proof.

\begin{theo}\label{therg1}
\begin{fleuvea}
\item A probability measure $\mu \in \mathcal{P}_{\text{inv}}(\Omega,
\mathcal{F})$ is extreme in $\mathcal{P}_{\text{inv}}(\Omega,
\mathcal{F})$ if and only if $\mu$ is ergodic.
\item Let $\displaystyle{\mu \in \mathcal{P}_{\text{inv}}(\Omega,
\mathcal{F}) \text{ and } \nu \in \mathcal{P}(\Omega, \mathcal{F})
\text{ such that } \nu \ll \mu,}$ then\\ $\displaystyle{ \nu \in
\mathcal{P}_{\text{inv}}(\Omega, \mathcal{F}) \text{ if and only if }
\exists \, h \geq 0, \; \mathcal{I} \text{-measurable}: \nu=h \mu}.$
\end{fleuvea}
\end{theo}

\begin{lemm}\label{lemerg1}
Let $\mu \in \mathcal{P}_{\text{inv}}(\Omega, \mathcal{F})$, then
$\mathcal{I} \subset \mathcal{F}_{-\infty} \; \mu$-a.s. More
precisely, for each $A \in \mathcal{I}$ there exists $B \in
\mathcal{F}_{-\infty}$ such that $\mu(A \Delta B) = 0$.
\end{lemm}
\textbf{Proof}\\
Let $A \in \mathcal{I}$ and $\left(B_{n}\right)_{n \geq 1}$ be a
sequence of cylinder sets such that $\mu(A \Delta B_{n}) \leq 2^{-n}$
for all $n \geq 1$. Since $\mu \in
\mathcal{P}_{\text{inv}}(\Omega,\mathcal{F})$, we have that
\[
\mu\left(A \Delta \tau^{i} B_{n}\right) \;=\;
\mu\left(\tau^{i} A \Delta \tau^{i} B_{n}\right) \;=\;
\mu\left(A \Delta B_{n}\right) \;\leq\; 2^{-n}
\]
for each $i \in\mathbb{N}$ ($\tau^{i}$ is the $i$th-iterate of
$\tau$).  Consider $\Lambda_{n} \uparrow \lat$ such that $B_{n} \in
\mathcal{F}_{\Lambda_{n}}$. For each $n \geq 1$ we choose $i(n) \geq
0$ such that $\Lambda_{n} \cap (\Lambda_{n}-i(n)) = \emptyset$. Each
set $C_{n} \triangleq \tau^{i(n)} B_{n}$ belongs to
$\mathcal{F}_{\left(\Lambda_{n}\right)_{-}}$ and satisfies $\mu\left(A
\Delta C_{n} \right) \leq 2^{-n}$. Therefore, the set $C \triangleq
\bigcap_{m \geq 1} \bigcup_{n \geq m} C_{n}$ belongs to
$\mathcal{F}_{-\infty}$ and satisfies
\[
\mu(A \Delta C) \;\leq\; \mu\left(\bigcap_{m
\geq 1} \bigcup_{n \geq m} A \Delta C_{n}\right) \;\leq\; \lim_{m
\rightarrow +\infty}\sum_{n \geq m} 2^{-n} \;=\; 0\;. \quad \Box
\]
\smallskip

\noindent\textbf{Proof of Theorem \protect\ref{th4}}\\
\noindent{\bf (a)} Let us consider the probability kernel $T$ on
$\mathcal{F} \times \Omega$ defined by
\[T(A \mid \omega) \;=\; \ind{A}(\tau \omega)\]
for every $A \in \mathcal{F}$ and every $\omega \in \Omega$.
\smallskip

To prove necessity we introduce
\[
\mathcal{K}(\mu) \;\triangleq\; \left(\bigcap_{\Lambda \in
  \mathcal{S}_{b}} \mathcal{I}_{f_{\Lambda}}^{\mathcal{F}_{\leq
    m_\Lambda}}(\mu) \right) \bigcap
\mathcal{I}_{T}^{\mathcal{F}}(\mu)\;.
\]
By \reff{lemext3} and Lemma \ref{lemerg1}, $\mathcal{K}(\mu)$ is the
$\mu$-completion of $\mathcal{I}$.  Therefore Corollary \ref{corext1}
implies that each $\mu$ extreme in $\mathcal{G}_{\text{inv}}(f)$ is
trivial on $\mathcal{I}$.

For the sufficiency, suppose that $\mu$ is trivial on $\mathcal{I}$
and consider a decomposition $\mu = s \, \nu + (1-s) \, \nu'$ with
$0<s<1$ and $\nu, \, \nu' \in \mathcal{G}_{\text{inv}}(f)$. Then there
exist $\mathcal{F}$-measurable $h, \, h' \geq 0$ such that $\nu = h
\mu$ and $\nu' = h' \mu $. Since $\mu, \, \nu, \, \nu' \in
\mathcal{P}_{\text{inv}}(\Omega, \mathcal{F})$, Proposition
\ref{prop1} applied to $\mathcal{I}_{T}^{\mathcal{F}}(\mu)$ implies
that $h, \, h'$ are measurable with respect to the $\mu$-completion of
$\mathcal{I}$. Hence the triviality of $\mu$ on $\mathcal{I}$ assure
that $h=h'=1 \; \mu$-a.s. Thus $\mu = \nu = \nu'$.
\smallskip

\noindent{\bf (b)}. Theorem \ref{therg1} (b) implies that there exists
$h \geq 0$, $\mathcal{I}$-measurable such that $\nu = h \mu$. By Lemma
\ref{lemerg1} $h$ is $\mathcal{F}_{-\infty}$-measurable, so Theorem
\ref{th2} b) implies that $\nu \in \mathcal{G}(f)$. Therefore $\nu \in
\mathcal{G}_{\text{inv}}(f)$.
\smallskip

\noindent {\bf (c)} It is an immediate consequence of (b). $\quad
\Box$

\section{Proofs on uniqueness}
\subsection{One-sided boundary-uniformity}
\begin{lemm}\label{lem2}
If uniqueness condition \reff{lis91} is satisfied, then $\nu \geq c \;
\mu, \; \forall \, \mu, \nu \in \mathcal{G}(f).$
\end{lemm}
\textbf{Proof}\\
Let $A$ be a cylinder set and $n$ an integer such that \reff{lis91}
holds. If $\mu$ and $\nu$ are consistent with $f$,
\begin{eqnarray*}
\nu(A) &=& \iint f_{[-n,m]}(A \mid \xi)\, \mu(d \eta)\, \nu(d \xi)\\
 & \geq & c \iint f_{[-n,m]}(A \mid \eta)\, \mu(d \eta)\, \nu(d \xi)\\
& = &c \, \mu(A). \quad \Box
\end{eqnarray*}
\smallskip

\noindent\textbf{Proof of Theorem \protect \ref{th6}}\\ We shall prove that
every element of $\mathcal{G}(f)$ is extreme. Let $\mu \in
\mathcal{G}(f)$ and $B \in \mathcal{F}_{- \infty}$ such that
$\mu(B)>0$. Define
\[
\nu \;\triangleq\; \mu( \cdot \mid B) \;=\;
\frac{\ind{B}}{\mu(B)} \, \mu\;.
\]
By Theorem \ref{th2} (c), $\nu \in \mathcal{G}(f)$. By the preceding
lemma $0 = \nu(B^{c}) \geq c \,\mu(B^{c})$, so
$\mu(B)=1. \quad \Box$
\medskip

\noindent\textbf{Proof of Proposition \protect \ref{cuc3}}\\
Call $m(f)$ the infimum \reff{eq:unif} and $V(f)$ the supremum
\reff{cuc5}.  Through an elementary logarithmic inequality we have
that for each $i,j\in\lat$ with $i>j$ and each
$\xi,\eta\in\Omega_{\le i}$ with $\xi_{j}^{i} = \eta_{j}^{i}$,
\equa{\label{eq:cuc.1.r}
\frac{f_{\{i\}}\bigl(\xi_i \bigm| \xi_{-\infty}^{i-1}\bigr)}{
f_{\{i\}}\bigl(\eta_i\bigm|\eta_{-\infty}^{i-1}\bigr)}
\;\ge\; \exp\Bigl(-\frac{\var{j}(f_{\{i\}})}{m(f)}\Bigr)\;.
}
Applying the factorization \reff{eq:facr} we conclude that for each
$n,m\in\lat$ with $n<m$ and each
$\xi,\eta\in\Omega_{\le m}$ with $\xi_{n}^{m} = \eta_{n}^{m}$,
\equa{\label{eq:cuc.2.r}
\frac{f_{[n,m]}\bigl(\xi_n^m\bigm|\xi_{-\infty}^{n-1}\bigr)}{
f_{[n,m]}\bigl(\eta_n^m\bigm|\eta_{-\infty}^{n-1}\bigr)}
\;\ge\; e^{-V(f)/m(f)}\;.\quad\Box
}

\subsection{Dobrushin uniqueness}
The following bound is the basic tool of the theory.

\begin{lemm}[Multisite dusting lemma]\label{duslem2}
Let $V \in \mathcal{S}_{b}$, $f_V$ a probability kernel on
$\mathcal{F}_{\le m_V}\times \Omega$ and
$\alpha^{V}$ is a $d$-sensitivity estimator for $f$.  Then,
\equa{\label{dus23}
\delta_{j}^{d}\left( f_{V}h\right) \quad
\begin{cases} =\;0 & \text{ if } j \in V\\[5pt]
\displaystyle{\leq\; \delta_{j}^{d}(h)+\sum_{k \in
    V}\delta_{k}^{d}(h) \, \alpha_{kj}^{V}} & \text{ if }
    j \in V_{-}\;,
\end{cases}
}
for every continuous function $h$ on $V \cup V_{-}.$
\end{lemm}
\begin{rema} The name of the lemma comes from a picturesque
  interpretation due to Michael Aizenman reported in Simon
  (1993\nocite{sim93}): If the oscillations are interpreted as ``dust"
  and the averages $f_{V}$ as applications of a (multisite) ``duster",
  the lemma says that no dust remains in $V$ after dusting the sites
  there [first line of \reff{dus23}], but the dust has been spread
  over the remaining sites [second line of \reff{dus23}].  The
  estimators give the fraction blown from site to site.  In this
  picture, Dobrushin condition \reff{gen30} means that some dust stays
  in the duster, a fact that allows for an eventual total cleaning.
\end{rema}

\noindent\textbf{Proof}\\
The first line in \reff{dus23} just
expresses the fact that the average $f_{V}h$ is
$\mathcal{F}_{V_-}$-measurable. The second line shows two
contributions: The first one due to the direct dependence of $h$ on
the configuration at the site $j$, and the second to the sensibility
of the $f_V$-averages to the configuration on the past instant $j$.
To separate both contributions we introduce a family of auxiliary
functions $h_{V, \omega}\left(\sigma_{V}\right) \triangleq
h\left(\omega_{V_{-}} \, \sigma_{V}\right)$ for each $\omega \in
\Omega$ (``freezing'' at $\omega$).  For $j \in V_-$ and $\xi, \, \eta
\in \Omega_{V_-}$ such that $ \xi \stackrel{\neq j}{=} \eta$, we have
\eqna{
\lefteqn{\Bigl| f_{V}(h \mid \xi) - f_{V}(h \mid \eta)\Bigr|}
\nonumber\\
&& \leq \; \Bigl|\stackrel{\circ}{f_{V}} \left(h_{V, \, \xi} -
h_{V, \, \eta}\mid \xi\right)\Bigr| +
\Bigl|\stackrel{\circ}{f_{V}} \left( h_{V, \,
\eta} \mid \xi\right) - \stackrel{\circ}{f_{V}}
\left(h_{V, \, \eta} \mid \eta \right)\Bigr|\;.
}
If we divide throughout by $d(\xi_j,\eta_j)$ and use the estimator
bound \reff{dus19} we obtain, upon taking the necessary suprema, the
second line in \reff{dus23}. $\quad \Box$\\

We now fix a partition $\mathcal{P}$ of $\lat$ into finite intervals and denote, for each $\Lambda \subset
\mathcal{S}_{b}$,
\[
\Lambda^*\;=\;\bigcup\,\bigl\{ V\in\mathcal{P} : \Lambda\cap V\neq\emptyset \bigr\}\;.
\]
Let $n(\Lambda)$ denote the number of elements of $\mathcal{P}$ forming $\Lambda^*$.

\begin{prop}\label{genth1}
Consider a LIS $f$ and $d$-sensitivity estimators $\alpha^{V}$ for
$f_{V}$ for each $V \in \mathcal{P}$.
\begin{fleuvei}
\item For every $j \in \Lambda_{-}^{*}$ and $h
\in \mathcal{B}_d(\Lambda^{*} \cup \Lambda_{-}^{*})$,
 \equa{\label{gen13}
\delta_{j}^{d}\left( f_{\Lambda^{*}}h \right) \;\leq\;
 \delta_{j}^{d}(h) + \sum_{k \in \Lambda^{*}} \delta_{k}^{d}(h)
\left[ \sum_{l=1}^{n(\Lambda)} \left( P_{\Lambda^{*}} \alpha
  \right)^{l}\right]_{kj}\;.
}
\item If Dobrushin condition \reff{gen30} is satisfied, then for every
  $j \in \Lambda_{-}^{*}$ and $h \in \mathcal{B}_d(\Lambda^{*})$,
\equa{\label{gen15} \delta_{j}^{d}\left( f_{\Lambda^{*}}h \right)
\;\leq\; \sum_{k \in \Lambda^{*}} \delta_{k}^{d}(h)
\left[ \frac{P_{\Lambda^{*}} \alpha }{1 - P_{\Lambda^{*}}\alpha}
  \right]_{kj}\;.
}
\end{fleuvei}
\end{prop}
\textbf{Proof}\\
We only need to prove \reff{gen13}.  Inequality \reff{gen15} is then
obtained by bounding the sum in the RHS of \reff{gen13} by the limit
$n(\Lambda)\to\infty$, which is finite under Dobrushin condition.
\smallskip

We proceed by induction on $n(\Lambda)$. The case $n(\Lambda)=1$ is just the multisite dusting lemma.  Suppose the
inequality valid for all $\Lambda$ with $n(\Lambda)=n$.  Consider $\Delta$ such that $\Delta^{*}=
\bigcup_{i=1}^{n+1}V_{i}$, where the $V_{i} \in \mathcal{P}$, $i=1, \ldots, n+1$ are labeled so that
$m_{V_i}=l_{V_{i+1}-1}$.  Denote $\Lambda^{*}= \bigcup_{i=1}^{n}V_{i}$.  Let $j \in \Delta_{-}^{*}$ and $h \in
\mathcal{B}_d(\Delta^{*}\cup \Delta_{-}^{*})$.  By the factorization property \reff{eq:facr2} of the LIS,
$\delta_{j}^{d}\left(f_{\Delta^{*}}h\right) = \delta_{j}^{d}\left(f_{\Lambda^{*}}f_{V_{n+1}}h\right)$.  Therefore, by
the inductive hypothesis,
\[
\delta_{j}^{d}\left(f_{\Delta^{*}}h\right) \;\leq\;
\delta_{j}^{d}(f_{V_{n+1}}h) + \sum_{k \in \Lambda^{*}} \delta_{k}^{d}
\left( f_{V_{n+1}}h \right)\left[ \sum_{l=1}^{n}\left(P_{\Lambda^{*}}
  \alpha \right)^{l} \right]_{kj}\;,
\]
and the multisite dusting lemma \ref{duslem2} yields
\begin{eqnarray*}
\lefteqn{\delta_{j}^{d}\left(f_{\Delta^{*}}h\right)  \;\leq\;
\delta_{j}^{d}(h) + \sum_{m \in V_{n+1}}\delta_{m}^{d}(h)
\left[P_{V_{n+1}}\alpha \right]_{mj}}\\
&& + \;\sum_{k \in \Lambda^{*}} \Biggl(\delta_{k}^{d}(h) + \sum_{m \in
  V_{n+1}}\delta_{m}^{d}(h)
\left[P_{V_{n+1}}\alpha \right]_{mk} \Biggr) \,
\left[\sum_{l=1}^{n}\left(P_{\Lambda^{*}} \alpha
\right)^{l}\right]_{kj}\;.
\end{eqnarray*}
We now observe that, given the restrictions in the sites being summed
over, we can replace in the RHS $P_{\Lambda^*}$ and $P_{V_{n+1}}$ by
  $P_{\Delta^*}$.  Furthermore, for $m\in V_{n+1}$, $l\in\mathbb{N}$,
\[
\sum_{i=1}^{n} \sum_{k \in V_{i}}
\left[P_{\Delta^{*}} \alpha \right]_{mk}
\left[ \left( P_{\Delta^{*}} \alpha \right)^{l}\right]_{kj} \;=\;
\left[ \left( P_{\Delta^{*}} \alpha \right)^{l+1}\right]_{mj}\;.
\]
The last two displays imply that
\[
\delta_{j}^{d}\left(f_{\Delta^{*}}h\right) \;\leq\;
\delta_{j}^{d}(h) + \sum_{k \in \Delta^{*}} \delta_{k}^{d}(h) \left[
  \sum_{l=1}^{n+1}\left( P_{\Delta_{*}} \alpha \right)^{l}
  \right]_{kj}\;. \quad \Box
\]

\noindent\textbf{Proof of Theorem \protect \ref{gencor1}}\\ Let us
label the elements of the partition so that
$\mathcal{P}=\{V_i:i\in\mathbb{Z}\}$ and $m_{V_i}=l_{V_{i+1}-1}$,
$i\in\mathbb{Z}$. Let us denote $V_{m-i}^{n}=\bigcup_{l=m-i}^{n}V_{l}$
for every integer $n,m,i$ with $m-i\le n$. Let $\mu, \, \nu \in
\mathcal{G}(f)$ and consider a local function $h$ of $d$-bounded
variations.  Pick $m,n \in \mathbb{Z}$ such that $h \in
\mathcal{B}_d(V_m^n)$.  The consistency of both $\mu$ an $\nu$ with
$f_{V_{m-i}^{n}}$, for an integer $i>0$, imply
\[
\Bigl|\nu(h) - \mu(h)\Bigr| \;\leq\; \iint
\left|f_{V_{m-i}^{n}}\left(h \mid \xi \right) -
f_{V_{m-i}^{n}}\left(h \mid \eta \right) \right| \, \nu(d \xi) \,
\mu(d \eta)\;.
\]
Therefore, by the continuity of $f$ and \reff{gen15},
\begin{eqnarray*}
\Bigl|\nu(h) - \mu(h)\Bigr| &\leq& \sum_{j \in
  \left(V_{m-i}\right)_{-}} \delta_{j}^{d}\left(f_{V_{m-i}^{n}} h
\right) \, \iint d(\xi_{j},\eta_{j}) \, \nu(d \xi) \, \mu(d \eta)\\
& \leq & D\,\sum_{k \in \Lambda} \, \delta_{k}^{d}(h) \sum_{j \in
  \left(V_{m-i}\right)_{-}}
\left[ \frac{P_{\Lambda} \alpha}{1- P_{\Lambda} \alpha} \right]_{kj}\;.
\end{eqnarray*}
Under condition \reff{gen30} the series on the RHS is summable, hence
the bound converges to zero as $i\to\infty$. $\quad \Box$

\section{Proofs on loss of memory and mixing} \label{s.plmm}
\textbf{Proof of Theorem \protect \ref{memth1} and Proposition
  \protect\ref{p.ff1}}\\
Part (i) of Theorem \ref{memth1} is just \reff{gen13}.  The triangular
property of $F$ implies that for each $i \in \Lambda^{*}$,
\[
\left[ \left(P_{\Lambda^{*}} \alpha\right)^{2}\right]_{kj} e^{F(k,j)}
\;=\; \sum_{i \in \Lambda^{*}} \alpha_{ki} \,
\alpha_{ij} \, e^{F(k,j)} \;\leq\; \sum_{i \in \Lambda^{*}}
\alpha_{ki} \, e^{F(k,i)} \, \alpha_{ij} \, e^{F(i,j)}\;.
\]
Therefore,
\[
\left[\left(P_{\Lambda^{*}} \alpha\right)^{2}\right]_{kj} e^{F(k,j)}
\;\le\; \sum_{j \in \lat } \left[ \left(P_{\Lambda^{*}}
  \alpha\right)^{2}\right]_{kj} \, e^{F(k,j)} \;\leq \;
\gamma_{\Lambda^{*} }^{2}\;.
\]
Proceeding inductively we obtain
\equa{\label{eq:pln}
\left[\left(P_{\Lambda^{*}} \alpha\right)^{n}\right]_{kj} \;\leq\;
\gamma_{\Lambda^{*}}^{n} \; e^{-F(k,j)}
}
for every natural $n$. This yields \reff{mem21} upon summation over
$n$. Combining \reff{mem21} with \reff{gen15}, we obtain
\reff{mem15}. $\quad \Box$
\bigskip

\noindent\textbf{Proof of Theorem \protect \ref{comth1}}\\
Fix $\Lambda \in \mathcal{S}_{b}$ and $h \in \mathcal{B}_d(\Lambda)$.
Using the consistency of $\mu$ and $\widetilde{\mu}$ respectively with
$f$ and $\widetilde{f}$, we have that, for each $n\in\mathbb{N}$,
\eqna{\label{com9}
\Bigl| \mu(h)- \widetilde{\mu}(h)\Bigr| & \leq&
 \Bigl|\mu\left(f_{[m_\Lambda-n,m_\Lambda]}h\right)-
\widetilde{\mu}\left( f_{[m_\Lambda-n,m_\Lambda]}h\right) \Bigr|
\nonumber\\
&&\;+\;
\Bigl|\widetilde{\mu}\left(f_{[m_\Lambda-n,m_\Lambda]}h \right)
- \widetilde{\mu}\left( \widetilde{f}_{[m_\Lambda-n, m_\Lambda]}
h\right)\Bigr|\;.
}
We estimate separately each term on the right as $n$ tends to
infinity.  The compactness of $\Omega$ implies that $
f_{[m_\Lambda-n,m_\Lambda]}(h\,|\,\omega) \to\mu(h)$ for each
$\omega\in\Omega$ as $n\to\infty$ (see Remark \ref{rem:r1}).
Therefore, by dominated convergence ($h$ is continuous, hence bounded)
\equa{\label{com11}
\Bigl|\mu\left(f_{[m_\Lambda-n,m_\Lambda]}h\right)-
\widetilde{\mu}\left( f_{[m_\Lambda-n,m_\Lambda]}h\right) \Bigr|
\;\tend{n}{\infty }{}\; 0\;.
}
\smallskip

To bound the last term in \reff{com9} we telescope using the
factorization property \reff{eq:facr2} for LIS:
\eqna{\label{com15}
\lefteqn{
\Bigl|\widetilde{\mu}\bigl(f_{[m_\Lambda-n,m_\Lambda]}h \bigr)
- \widetilde{\mu}\bigl(\widetilde{f}_{[m_\Lambda-n,m_\Lambda]}
h\bigr) \Bigr| \;\leq\;
\Bigl|\mu\bigl(f_{\{m_\Lambda\}} h \bigr)
- \widetilde{\mu}\bigl(\widetilde{f}_{\{m_\Lambda\}} h \bigr)\Bigr|
\qquad}
\nonumber\\
&& \qquad \qquad\qquad
{} + \sum_{k=m_\Lambda-n}^{m_\Lambda-1}\Bigl|\widetilde{\mu}
\bigl(f_{[k,m_\Lambda]}h \bigr)-
\widetilde{\mu}\bigl(\widetilde{f }_{\{k\}}f_{[k+1,m_\Lambda]}
h\bigr)\Bigr| \;.
}
The definition
    \reff{vkr}/\reff{vkr2} of the VKR distance, implies that
\[
\Bigl| (f_{k}\,g)(\omega) - (\widetilde{f}_{k}\,g)(\omega) \Bigr|
\; \leq\; \delta_{k}^{d}(g) \; \Bigl\| \stackrel{\circ}{f_{\{k\}}}
\left( \, \cdot \mid \omega \right) -
\stackrel{\circ}{ \widetilde{f}_{\{k\}}}
\left( \, \cdot \mid \omega \right)\Bigr\|_{d}\;,
\]
for all $k \in \lat$, $\omega \in \Omega_{-\infty}^{k-1}$ and $g\in\mathcal{B}_d(]-\infty,k])$.  Hypothesis \reff{com5}
implies
\equa{\label{com17}
\Bigl| \widetilde{\mu} \bigl( f_{\{k\}}\,g - \widetilde{f}_{\{k\}}\,g
\bigr) \Bigr| \;\leq \;
\widetilde{\mu} \left( b_{k} \right) \delta_{k}^{d} (g)\;.
}
Combining \reff{com15} and \reff{com17} we obtain
\eqna{\label{eq:nti}
\lefteqn{
\Bigl|\widetilde{\mu}\bigl(f_{[m_\Lambda-n,m_\Lambda]}h \bigr)
- \widetilde{\mu}\bigr( \widetilde{f}_{[m_\Lambda-n,m_\Lambda]} h
\bigr)\Bigr|\; = \;}\nonumber\\
&&\qquad\qquad
\sum_{k=m_\Lambda-n}^{m_\Lambda-1}\widetilde{\mu}\left( b_{k}\right)
\,\delta_{k}^{d}\bigl( f_{[k+1,m_\Lambda]}\bigr) \;+\;
\widetilde{\mu}\left(b_{i}\right) \, \delta_{i}^{d}(h)\;.
}
To obtain \reff{com7} we insert this bound in \reff{com9}, let $n$
 tend to infinity and use \reff{com11}. $\quad \Box$
\bigskip

\noindent\textbf{Proof of Theorem \protect \ref{corth1}}\\
Fix $\Lambda, \Delta \in \mathcal{S}_{b}$ with $m_\Delta<l_\Lambda$,
$h_{1} \in \mathcal{B}_d(\Lambda)$ and $h_{2} \in
\mathcal{B}_d(\Delta)$. Without
loss, we can suppose that $h_{2} \geq 0, \; h_{2} \not\equiv 0$ and
$\mu(h_{2})=1$ since both sides of \reff{cor5} are invariant under
adding a constant to $h_{2}$ and both multiply in the same way if
$h_{2}$ is multiplied by a positive constant.   We then can write
\equa{\label{eq:apu10}
\text{Cor}_{\mu}\left(h_{1},h_{2}\right) \;=\; \Bigl|
\nu(h_{1}) - \mu(h_{1})\Bigr|
}
where $\nu$ is the probability measure defined by
\equa{\label{cor7}
\nu \;=\; h_2\,\mu\;.
}

\noindent\emph{1st stage:} We construct a LIS $\widetilde f$ for $\nu$
on $]-\infty,m_{\Lambda}]$. For every $k \in ]-\infty,m_{\Lambda}]$,
let us define
\equa{\label{cor11}
\widetilde{f}_{k} \;=\; g_{k}\, f_{\{k\}}
}
with
\equa{\label{cor15}
g_{k}\; =\; \begin{cases}
1 & \quad \text{if } k \in [m_{\Delta}+1,m_{\Lambda}]\\[5pt]
{\displaystyle \frac{f_{[k+1,m_\Delta]}\,
\bigl(h_{2}\bigm|\cdot\,\bigr)}{f_{[k,m_\Delta]}
\bigl( h_{2}\bigm|\cdot\,\bigr)} }& \quad
\text{if } k \in ]-\infty,m_{\Delta}]\;.
\end{cases}
}

The function $g_{k}$ is well defined because $f_{[k,m_\Delta]}h_{2}
\neq 0$ for every $k \in ]-\infty,i]$.  Indeed the existence of $k$
such that $f_{[k,q]}h_{2}=0$ would imply, by consistency, that
$\mu(h_{2})=0$.  This contradicts the fact that $\mu(h_{2})=1.$ It is
clear that the kernels $\widetilde{f}_{k}$ satisfy the hypotheses of
Theorem \reff{listh1}, hence they uniquely define a LIS
$\widetilde{f}$ on $]-\infty,m_{\Lambda}]$.  The same theorem shows
that the consistency of $\nu$ with each $\widetilde{f}_{k}$, $k \in
]-\infty,m_{\Lambda}]$ is all that has to be checked in order to prove
that $\nu$ is consistent with $\widetilde{f}$.
\smallskip

If $k \in [m_{\Delta}+1,m_{\Lambda}]$, this consistency is a
consequence of the following sequence of identities, valid for every
$h \in \mathcal{F}_{\leq k}$:
\equa{\label{eq:con1.r}
\nu \bigl( \widetilde{f}_{k}(h)\bigr)\; =\;
\mu(h_{2}\,f_{\{k\}}(h)\bigr) \;=\;
\mu(f_{\{k\}}(h_{2}\,h)\bigr) \;=\; \mu(h_2\,h)\;=\;\nu(h)\;.
}
The third inequality is due to the $\mathcal{F}_{\le k-1}$
measurability of $h_2$ and the fourth one to consistency.

For $k \in ]- \infty,m_{\Delta}]$ we observe that for
$h \in \mathcal{F}_{\leq k}$,
\[
\nu \bigl( \widetilde{f}_{k}(h)\bigr)\; =\;
\mu(h_{2}\,f_{\{k\}}(g_k\,h)\bigr) \;=\;
\mu\Bigl(f_{[k,m_\Delta]}\bigl[h_{2}\,f_{\{k\}}(g_k\,h)\bigr]\Bigr)\;,
\]
the last inequality being a consequence of the consistency of $\mu$
with $f$.  Upon inserting the definition of $g_k$ [second line in
  \reff{cor15}] we see that there is a term $f_{[k,m_\Delta]}$ in the
  denominator that can be pulled to the left because of its
  $\mathcal{F}_{\le k-1}$-measurability.  This produces a cancellation
  with an analogous term in the numerator.  We thus obtain
\eqna{\label{eq:con2.r}
\nu \bigl( \widetilde{f}_{k}(h)\bigr)& =&
\mu\Bigl(f_{\{k\}}\bigl[h\,f_{[k,m_\Delta]}(h_{2})\bigr]\Bigr) \;=\;
\mu\Bigl(f_{\{k\}}\bigl[f_{[k,m_\Delta]}(h_{2}\,h)\bigr]\Bigr)
\nonumber\\
&=& \mu(h_2\,h) \;=\; \nu(h)\;.
}
The third inequality is due to the $\mathcal{F}_{\le k}$-measurability
of $h$ and the fourth one to the consistency of $\mu$ with $f$.
Identities \reff{eq:con1.r} and \reff{eq:con2.r} prove that $\nu$ is
consistent with $\widetilde{f}$ on $]-\infty,m_{\Lambda}]$.
\smallskip

\noindent\emph{2nd stage:} For every $k
\in \Lambda \cup \Lambda_{-}$ and $\omega \in \Omega_{-\infty}^{k-1}$,
we construct $b_{k}(\omega)$ such that
\equa{\label{eq:apu}
\Bigl\| \stackrel{\circ}{f_{k}}(\, \cdot \mid \omega) -
\stackrel{\circ}{\widetilde{f}_{k}}( \, \cdot \mid \omega)
\Bigr\|_{d} \;\leq\; b_{k}(\omega)\;.
}
For starters, we can take
\equa{\label{cor17} b_{k}=0 \quad \forall \; k \in
[m_{\Delta}+1,m_{\Lambda}]\;,
}
because $\displaystyle{ \; \stackrel{\circ}{f_{k}}(\, \cdot \mid
  \omega) = \stackrel{\circ}{\widetilde{f}_{k}}(\, \cdot \mid
  \omega)}$, for $k \in [m_{\Delta}+1,m_{\Lambda}]$ and $ \omega \in
\Omega_{-\infty}^{k-1},$.

We fix $k \in \Delta \cup \Delta_{-}$ and $\omega \in
\Omega_{-\infty}^{k-1}$ and consider the set
$\Omega^\omega_k=\{\omega_k\in\Omega_{\{k\}} :
\omega_{-\infty}^k\in\Omega_{-\infty}^k\}$ with the restricted
topology and Borel $\sigma$-algebra.  To abbreviate the notation we
introduce the function $u:\Omega^\omega_k\to\mathbb{R}$ defined by
\equa{\label{cor19}
u(x) \;\triangleq\; g_{k}\left(\omega_{-\infty}^{k-1} \, x \right)
\;=\; \frac{f_{[k+1,m_{\Delta}]}\bigl(h_{2}\bigm|
\omega_{-\infty}^{k-1} \, x \bigr)}{
f_{[k,m_{\Delta}]}\bigl(h_{2}\bigm|\omega\bigr)}
}
and the measure
\equa{\label{cor21}
\alpha \;\triangleq\; \stackrel{\circ}{f_{k}}
\left( \, \cdot \mid \omega\right)
}
on $\Omega^\omega_k$.  Notice that
\equa{\label{eq:uau}
\stackrel{\circ}{\widetilde{f}_{k}}(\, \cdot \mid \omega)\, -
\stackrel{\circ}{f_{k}}(\, \cdot \mid \omega)\;=\;
u\, \alpha - \alpha\;.
}
We also denote, for each $\mathcal{F}_{\{k\}}$-measurable function
$h$,
\[
m_{h} \;\triangleq\; \sup_{x \neq y} \frac{h(x)+h(y)}{2d(x,y)}\;,
\]
and observe that
\equa{\label{cor23}
\bigl\|h-m_{h}D\bigr\|_{\infty} \;\leq\; \frac{D}{2} \, \delta_k^d(h)
}

We affirm that
\equa{\label{eq:apu2}
\Bigl\| \stackrel{\circ}{f_{k}}(\, \cdot \mid \omega) -
\stackrel{\circ}{\widetilde{f}_{k}}( \, \cdot \mid \omega)
\Bigr\|_{d} \;\leq\; \frac{D}{2}\, \alpha(|u-1|)\;.
}
Indeed, for $h\in\mathcal{B}_d(\{k\})$ with $\delta_k^d(h)\le
1$ we have
\[
\Bigl|u\, \alpha(h)-\alpha(h)\Bigr| \;=\;
\Bigl|\alpha\bigl[(u-1)(h-m_{h}D)\bigr]\Bigr| \;\leq\;
 \alpha \left( \left| u-1 \right| \right)\, \bigl\| h-m_{h}D
\bigr\|_{\infty}\;.
\]
From this and \reff{cor23}, assertion \reff{eq:apu2} follows.

We now use Schwarz's inequality to bound
\[
\alpha(|u-1|) \;=\; \alpha\bigl(|u-\alpha(u)|\bigr) \;\le\;
\Bigl[\alpha \Bigl( \left(u - \alpha(u)\right)^{2} \Bigr)
\Bigr]^{\frac{1}{2}}\;,
\]
and since $\alpha(u)$ minimize $x \longmapsto \alpha\left( \left(u-x
\right)^{2}\right),$ we obtain
\equa{\label{eq:apu5}
\alpha(|u-1|) \;=\;
\Bigl[\alpha \Bigl( \left(u - m_uD\right)^{2} \Bigr)
\Bigr]^{\frac{1}{2}} \;\le\; \|u-m_uD\|_\infty\;.
}
The combination of \reff{eq:apu2} and \reff{eq:apu5} gives
\reff{eq:apu} with
\equa{\label{eq:apu6}
b_{k}(\omega) \;\triangleq\;
\frac{D^{2} \, \delta_{k}^{d}(u)}{4} \;=\;
\frac{D^{2} \, \delta_{k}^{d}\left(f_{[k+1,m_{\Delta}]} h_{2}
  \right)}{4\, f_{[k,m_{\Delta}]}h_{2}(\omega)}\;.
}
\smallskip

\noindent\emph{3rd stage:}  We estimate $\nu\left(b_{k}\right)$. From
\reff{eq:apu6}:
\[
\nu(b_{k}) \;=\; \frac{D^{2}}{4}\, \delta_{k}^{d}\left(
  f_{[k+1,m_{\Delta}]} h_{2}\right) \; \mu \left( \frac{h_{2}}{
    f_{[k,m_{\Delta}]} h_{2}} \right)\;.
\]
By consistency, $\mu=\mu\,f_{[k,m_{\Delta}]}$, hence the last factor
  is just 1.  From this and \reff{cor17} we conclude that
\equa{\label{cor25}
\nu(b_{k}) \;=\; \begin{cases}
0 &\mbox{ if } k \in [m_{\Delta}+1,m_{\Lambda}]\\[5pt]
{\displaystyle  \frac{D^{2}}{4} \, \delta_{k}^{d}\left(
f_{[k+1,m_{\Delta}]} h_{2}\right)} & \mbox{ if }
k \in \Delta \cup \Delta_{-}\;.
\end{cases}
}
In view of \reff{eq:apu10}, \reff{eq:apu} and \reff{cor25}
imply \reff{cor5} by Theorem \ref{comth1}. $\quad \Box$
\bigskip

\noindent\textbf{Proof of Corollary \protect \ref{corcor1}}\\
Part (i) follows from \reff{gen15} and \reff{cor5}, and part (ii) from
\reff{gen13} and \reff{cor29}. $\quad \Box$

\section*{Acknowledgements}  R.F.\ wishes to thank Michael Aizenman
for many enlightening discussions on Dobrushin criterion.  The authors
are thankful to J.~Steiff for pointing out a mistake in a previous
version of this paper.

%
\bibliographystyle{alpha}
%

%
\end{document}